\documentclass[10pt]{amsart}
\usepackage{amssymb,url,mathrsfs,euscript,a4wide}

\theoremstyle{remark}     
\newtheorem{rmk}{Remark}[section]\newtheorem*{rmk*}{Remark}
\newtheorem*{rmk1}{Remarks 1}\newtheorem*{rmk2}{2}
\newtheorem*{ex*}{Example}
\newtheorem*{exs*}{Examples}

\newtheorem*{acknowledgements}{Acknowledgements}

\theoremstyle{plain}      
\newtheorem{lemma}[rmk]{Lemma}
\newtheorem{proposition}[rmk]{Proposition}
\newtheorem{theorem}[rmk]{Theorem}
\newtheorem{corollary}[rmk]{Corollary}

\theoremstyle{definition} 
\newtheorem{definition}[rmk]{Definition}

\numberwithin{equation}{section}

\newcommand{\bdm}{\begin{displaymath}}
\newcommand{\edm}{\end{displaymath}}

\newcommand{\hodge}{{*}}
\newcommand{\lie}[1]{\mathfrak{#1}}
\newcommand{\Lie}[1]{\textsl{#1}}
\newcommand{\G}{G\sb 2}
\DeclareMathOperator{\SO}{\Lie{SO}}
\DeclareMathOperator{\SL}{\Lie{SL}}
\DeclareMathOperator{\GL}{\Lie{GL}}

\DeclareMathOperator{\SU}{\Lie{SU}}
\DeclareMathOperator{\U}{\Lie{U}}

\DeclareMathOperator{\so}{\lie{so}}

\newcommand{\g}{\lie{g}}
\renewcommand{\l}{\lie{l}}

\newcommand{\s}{\lie{s}}
\newcommand{\h}{\lie{h}}
\renewcommand{\k}{\lie{k}}
\newcommand{\p}{\lie{p}}

\DeclareMathOperator{\Hol}{\textsl{Hol}\,}
\DeclareMathOperator{\Ric}{\textsl{Ric}}

\newcommand{\iso}{\cong}

\newcommand{\lra}{\longleftrightarrow}

\newcommand{\tto}{\mapsto}
\newcommand{\lan}{\langle}
\newcommand{\ran}{\rangle}
\newcommand{\hook}{\lrcorner\,}
\newcommand{\tsum}{\tx\sum}

\newcommand{\sym}{\mathcal{S}}

\newcommand{\1}{e^1}\newcommand{\2}{e^2}\newcommand{\3}{e^3}
\newcommand{\4}{e^4}\newcommand{\5}{e^5}

\renewcommand{\leq}{\leqslant}

\newcommand{\e}[1]{e_{#1}}

\newcommand{\ad}{\textsl{ad}}

\newcommand{\tx}{\textstyle}
\newcommand{\q}{\quad}\newcommand{\qq}{\qquad}

\newcommand{\ph}{\hphantom}
\newcommand{\bb}{\mathbb}

\newcommand{\bsub}{\begin{subequations}}
          \newcommand{\esub}{\end{subequations}}
\newcommand{\be}{\begin{equation}}\newcommand{\ee}{\end{equation}}
\newcommand{\ba}{\begin{array}}\newcommand{\ea}{\end{array}}
\newcommand{\bg}{\begin{gathered}}\newcommand{\eg}{\end{gathered}}

\newcommand{\gh}{\mathfrak}

\newcommand{\R}{\mathbb{R}}\newcommand{\C}{\mathbb{C}}

\newcommand{\w}{\wedge}
\newcommand{\set}{\subseteq}\newcommand{\sset}{\subset}
\newcommand{\na}{\nabla}

\newcommand{\vs}{\vphantom{\int_W^M}}

\def \NI{\textsl{NI}}

\def\sideremark#1{\ifvmode\leavevmode\fi\vadjust{\vbox to0pt{\vss
     \hbox to 0pt{\hskip\hsize\hskip1em%
     \vbox{\hsize2cm\tiny\raggedright\pretolerance10000%
     \noindent #1\hfill}\hss}\vbox to8pt{\vfil}\vss}}}%

\def\polhk#1{\setbox0=\hbox{#1}{\ooalign{\hidewidth
       \lower1.5ex\hbox{`}\hidewidth\crcr\unhbox0}}}

\allowdisplaybreaks

\begin{document}

\title[]{Nearly integrable $\boldsymbol{\SO(3)}$ structures \\
on  $\boldsymbol{5}$-dimensional Lie groups}

\author{Anna Fino} \address[A.Fino]{Dipartimento
       di Matematica, Universit\`a di Torino, via Carlo Alberto 10, 10123
       Torino, Italy}\email{fino@dm.unito.it}
\author{Simon G.~Chiossi} \address[S.Chiossi]{Institut f\"ur Mathematik,
       Humboldt-Universit\"at zu Berlin, Unter den Linden 6,
       10099 Berlin, Germany}\email{sgc@math.hu-berlin.de}

\begin{abstract}
Recent work \cite{Bobienski-N:SO(3)} on  $5$-dimensional Riemannian 
manifolds with an $\SO(3)$ structure prompts
us to investigate which Lie groups admit such a geometry. The case
in which the $\SO(3)$ structure admits a compatible connection with
torsion is considered. This leads to a classification under special 
behaviour of the connection, 
which enables to recover all known
examples, plus others bearing torsion of pure type.
Suggestive relations with 
special structures in other dimensions are highlighted, with attention
to eight-dimensional $\SU(3)$ geometry.
\end{abstract}

\subjclass[2000]{Primary 53A40 -- Secondary 53C10, 53B15, 53C35, 53C25}
\maketitle


\section{Introduction}
\noindent
Given an oriented Riemannian manifold of dimension five $(M^5,g)$, an $\SO(3)$
structure is the reduction of the structure group of the frame bundle
to the Lie group $\SO(3)$ sitting inside $\SO(5)$. The inclusion we
choose is the one based on the
 \emph{irreducible} $5$-dimensional representation of $\SO(3)$,
determined by the decomposition
\[
\so(5)=\so(3)\oplus V, 
\]
where $V$ is
the unique irreducible $7$-dimensional representation of $\SO(3)$.
It is known that the homogeneous space $\SO(5)/\SO(3)$ has an 
$\SO(3)$-connection $\tilde\nabla$
whose torsion tensor $T$ is skew-symmetric (and unique), by
which one intends that
\begin{equation}
\label{eq:T3}
T(X,Y,Z)=g(\tilde\na_XY-\tilde\na_YX-[X,Y],Z)
\end{equation}
 is a three-form for any $X,Y,Z$ vector fields, for which refer to
\cite{Agricola-F:holonomy,Friedrich:non-integrable}. Indeed this was
  already noticed in \cite{Salamon:holonomy}, to the effect that the
  skew-symmetry of $T$ can be used to align $\SO(3)$ structures on five-manifolds 
with other kinds of torsion geometries.
In the general framework of geometric structures on Riemannian
$n$-manifolds, $G$-reductions are distinguished by
the irreducible components of the representation
$\R^n\otimes\so(n)/\g$, $\g$ being the Lie algebra of $G$.
The decomposition of this space depends on a
tensorial object, most of the times a differential form.
Well known is the archetypal description of almost Hermitian
structures \cite{Gray-H:16} in terms of the K\"ahler
form.
Other $G$-geometries have been discussed
using the same approach, see
\cite{Cabrera:specialAH} amongst others
for $G=\SU(n)$.

We review $\SO(3)$ geometry in section \ref{sec:so3}. In some sense it
is slightly different from more familiar
$G$-structures, for it is defined not by means
of a skew form - 
but rather a symmetric tensor, denoted ${\mathbb T}$.
 Apart from ${\mathbb T}$, one can not
expect to find other invariants for the representations of $SO(3)$
on $\R^5$, in particular no differential form. 
A rather diverse situation occurs
in higher dimensions, such as eight for
instance, where both totally symmetric and totally skew $3$-forms play a
relevant role \cite{Nurowski:distinguished-dimensions}, see section 
\ref{sec:8}.

In the special case of concern, the so-called \emph{nearly integrable}
structures, this tensor behaves - strikingly - like the
almost complex structure $J$ of a nearly K\"ahler manifold, whence the
similar name.

For a general $\SO(3)$-reduction the tensor product $\R^5\otimes {\so(3)}^\perp$
decomposes into the analogues of the Gray--Hervella classes, minding
 that here the complement ${\so(3)}^\perp$ is reducible (hence not $V$).
As the torsion $T$ of a nearly integrable $\SO(3)$
structure is uniquely defined and skew-symmetric, and the spaces of
 three- and two-forms on $M^5$ are Hodge-isomorphic, one decomposes 
 the latter under $\SO(3)$ in the sum of two irreducible modules
\[\Lambda^2M^5=\Lambda^2_3 \oplus \Lambda^2_7
\] 
of dimensions three, seven respectively.
When the components therein are trivial,
the characteristic
connection $\tilde \na$ is in fact the Riemannian one, and \cite{Bobienski-N:SO(3)}
proves that the 
simply-connected manifolds $M^5$ admitting
a torsion-free $\SO(3)$ structure are
$\R^5$, $\SU(3)/\SO(3)$ or $\SL(3,\R)/\SO(3)$.
It is not surprising that these models are symmetric spaces, as 
proven by Berger's holonomy
theorem. The paper of Bobie\'nski and Nurowski classifies
$5$-dimensional Lie groups with type $\Lambda^2_3$, and finds examples  
with closed torsion of the complementary type $\Lambda^2_7$. 

The present note
pertains to $5$-dimensional real connected Lie groups 
$L$ with invariant metric
$g$ and $\SO(3)$ structure having anti-symmetric torsion $T$.
More precisely, 
we characterise the Lie groups $(L, \bb{T},g)$ admitting a
splitting $\l=\h\oplus \p$ of the Lie algebra of $L$ defined
via an adapted frame, 
see section \ref{sec:symmetric pair}, in terms of the structure 
constants. 
This in turn yields a classification (theorem
\ref{thm:main}) assuming that the
$\SO(3)$ connection $\tilde \nabla$ satisfies 
\[
\begin{array}{l}
\tilde\nabla_X {\gh h} \subseteq {\gh h}, \quad \tilde\nabla_X {\gh p}
\subseteq {\gh p}, \quad \forall X \in {\gh h}\\
\tilde\nabla_Y {\gh h} \subseteq {\gh p},  \quad \tilde\nabla_Y {\gh p}
\subseteq {\gh h}, \quad \forall Y \in {\gh p}.
\end{array}
\]
This special behaviour is displayed by the Levi--Civita
connection $\na$ in all instances of \cite{Bobienski-N:SO(3)}, and
corresponds to demanding that the group $L$ act transitively
on a Riemannian symmetric surface.
We then prove that either the Levi-Civita connection fulfills the same
algebraic conditions as $\tilde \nabla$, or 
$\tilde \nabla$ is identically zero (theorem \ref{thm:non-symm-T=0}).
In the latter case the Lie group is essentially 
$\SO(3)\times\R^2$, and thus the unique instance up to isomorphisms.
 Our results are in line with the main reference \cite{Bobienski-N:SO(3)}
and in some sense attempt to complete that classification. 
All examples have $d\hodge T=0$, and it would be thus natural to ask
whether this were always the case. We
prove that a large class of Lie groups,
comprising those of \cite{Bobienski-N:SO(3)}, satisfies the equivalent requirement of
symmetry of the $\SO(3)$ connection's Ricci tensor.
  In every case the
torsion type,
whether $\Lambda^2_3, \Lambda^2_7$ or generic, is determined. 
In section \ref{sec:example} we prove 
the existence of non-strong 
structures of type $\Lambda^2_7$, by finding an explicit example. This  
is once again realised by a Lie group that acts transitively on a
$3$-dimensional symmetric space by way of 
$$
[\p, \p] \subseteq \p, \quad [\h, \p] \subseteq \h, \quad [\h, \h] \subseteq \p.
$$

Section \ref{sec:8} is devoted to understanding the geodesic equation $
\nabla_X X =
0$, for any $X$ in $\l$, which arises from `near
integrability' and has remarkable consequences. Using this
property it is possible to construct nearly integrable  geometry of
higher dimension, namely $\SU(3)$ structures
on the Riemannian products $L \times \R^3$ and $L\times \SO(3)$.
\vspace{0.4cm}
\begin{acknowledgements}
Thanks are due to Th.~Friedrich for presenting the problem to us in the first
place and many, many discussions and to P.~Nurowski
for sharing an early version of
\cite{Bobienski-N:SO(3)}, answering questions and his constant
interest.
 Our gratitude goes to S.~Garbiero, S.~Salamon and A.~Swann for suggesting
 improvements to the manuscript.
SC thanks the Department of Mathematics at the University of Torino
for kind hospitality in early 2006.
The work has been supported by {\sc gnsaga} of {\sc in}d{\sc am}, {\sc
  miur}, and the {\sc sfb} 647
``Space--Time--Matter'' of the {\sc dfg}.

\medbreak
\noindent 
The authors would like to dedicate this paper to the memory of Giulio
Minervini (University of Bari).
\end{acknowledgements}

\section{Irreducible $\SO(3)$ geometry}
\label{sec:so3}

\noindent
The vector space $\R^5$ is isomorphic to the set of real $3\times 3$
symmetric matrices
with no trace $\sym^2_0\R^3$; we fix the isomorphism as follows
\begin{equation}
       \label{eq:rep}
       X=(x_1,\ldots,x_5)\lra
       \begin{pmatrix}
         \dfrac{x_1}{\sqrt{3}}-x_4 & x_2 & x_3\\
         x_2 &\dfrac{ x_1}{\sqrt{3}} +x_4 & x_5\\
         x_3 & x_5 & -\dfrac{2}{\sqrt{3}}x_1
       \end{pmatrix},
\end{equation}
the square root only being a convenient factor.
The irreducible representation on $\R^5$ is given by
$$
\rho(h) X = h X h^{-1},\q h\in\SO(3),
$$
 where the vector $X$ is thought of as matrix of the above form.
An $\SO(3)$
structure on $(M,g)$ can be identified
\cite{Bobienski-N:SO(3)} with an element of ${\mathbb T}\in\bigotimes^3\R^5$
such that\medbreak

i)\q ${\mathbb T}=
\sum_{i,j,k=1}^5 t_{ijk}\,dx_i\otimes dx_j\otimes dx_k$ is
symmetric 
in all arguments,\smallbreak

ii)\q the endomorphisms $X\tto {\mathbb T}_X=X\hook {\mathbb T}$ are
trace-free and\smallbreak

iii)\q $({\mathbb T}_X)^2X=g(X,X)X$, for all $X\in\R^5$.\smallbreak

\noindent
The contraction $X\hook {\mathbb T}$ prescribes to fill the first argument
of ${\mathbb T}$, so ${\mathbb T}_X(\cdot,\cdot)={\mathbb
  T}(X,\cdot,\cdot)$. 
Notice by the way that ii) is expected, given
that $\sum_jt_{ijj}$ defines the components of an $\SO(3)$-invariant vector in $\R^5$,
so it must vanish for each $i=1\ldots 5$ because the representation is irreducible.

These properties together with the identification \eqref{eq:rep}
determine an adapted frame, i.e.~an orthomormal basis $\{e_1,e_2,e_3,e_4,e_5\}$ of tangent
vectors so that $X=\sum_{j = 1}^5 x_je_j$ and the metric automatically assumes
the canonical form $g=\sum_{i = 1}^5 (e^i)^2$.
Adopting the choice made in \cite{Bobienski-N:SO(3)}, the $\SO(3)$
structure is completely described by the tensor
\begin{equation} \label{eq:t_ijk}
{\mathbb T}(X,X,X)= x_1(3x_2^2+3x_4^2-x_1^2-\tfrac 32 x_3^2-\tfrac 32
x_5^2)+\tfrac {3\sqrt{3}}2 x_4(x_5^2-x_3^2)+ 3\sqrt{3}x_2x_3x_5,
\end{equation}
given by the determinant of $X$.
\bigbreak
The action of $\SO(3)$ on $\C^2$, like the one of $\SU(2)$,
 endows the symmetric tensor product 
$\sym^4\C^2\iso\C^5$ with a real structure, and 
$\sym^3(\sym^4\C^2)$ contains only one copy of $\C$,
generated in fact by ${\mathbb T}$. The recipe to
tackle a general $G$-structure prescribes first to determine the form(s)
with isotropy $G=\SO(3)$, the lack of which 
makes this whole matter quite complicated.

\begin{definition}
\cite{Bobienski-N:SO(3)}
An $\SO(3)$ structure is said \emph{nearly integrable} (\NI) when
\begin{equation} \label{eq:ni}
(\nabla_X{\mathbb T})(X,X,X)=0
\end{equation}
for all vector fields $X\in\R^5$.
\end{definition}

\noindent
This relation is -- at least
formally -- similar to that defining a nearly K\"ahler
structure on an even dimensional manifold: $(\na_X J)X=0$. Just as the
almost complex structure 
$J$ is a Killing form \cite{Semmelmann:Killing-forms} there,
\eqref{eq:ni} is saying that 
${\bb T}$ is a symmetric Killing $3$-tensor. 

The theoretical interest of such a structure lies in the fact that it
admits a uniquely defined \emph{characteristic} connection
\begin{equation}\label{eq:torsion}
\tilde\nabla=\nabla-\tfrac 12 T
\end{equation}
with torsion $T$ a three-form, in the sense of
\eqref{eq:T3}. We finally ought to remind that
the existence of a
nearly integrable $\SO(3)$ structure is actually the same
\cite{Bobienski-N:SO(3)} as having 
\emph{skew-symmetric torsion} $T$. 
Among the simplest examples of $G$-invariant metric connections with
anti-symmetric torsion one counts those of naturally
reductive spaces. These are homogeneous spaces with a reductive
decomposition $\h\oplus\p$ possessing a connection 
$\tilde\na=\na+[\cdot,\cdot]_\p$ and skew-symmetric
characteristic torsion tensor $-g([\cdot,\cdot]_\p,\cdot)$ by very definition. 
Because the latter and the induced
curvature are parallel, naturally
reductive spaces generalise Riemannian symmetric spaces.
For a $G$-structure though,
there is no Lie group acting transitively.
Besides, there are sound reasons to believe that 
Riemannian manifolds $M^n$ endowed with metric connections $\tilde\na$ 
with skew torsion are useful in string and supergravity theories, see 
\cite{Friedrich-I:skew, Agricola:srni} and references.
\bigbreak
The space $\Lambda^3\R^5$ is isomorphic via the Hodge operator
to that of two-forms,
and the latter decomposes under $\SO(3)$ into the direct sum of
the irreducible modules
\begin{equation}
\label{eq:torsion-types}
\Lambda^2_3 =\text{span}\{E_1,E_2,E_3\}\iso\so(3),\qq
\Lambda^2_7 = \bigl(\Lambda^2_3 \bigr)^\perp.
\end{equation}
The orthogonal complement 
is taken with respect to the pairing
\begin{equation}\label{eq:torsion-pairing}
\lan
\alpha,\beta\ran=\hodge(\hat {\mathbb T}(\alpha)\w\hodge\beta),
\end{equation}
    via the endomorphism
$\hat {\mathbb T}(e^i\w e^k)=4\sum_{j,l,m} t_{ijm}t_{klm}e^j\w e^l$ of
$\Lambda^2(\R^5)^*$
from \cite{Bobienski-N:SO(3)}. Wedge products of basis vectors will be
henceforth expressed by juxtaposed indexes. 

The forms $E_i$ are given in matrix form by
\bdm
E_1={\tiny\begin{pmatrix}
0&0&0&0&\sqrt{3}\\
0&0&1&0&0\\
0&-1&0&0&0\\
0&0&0&0&1\\
-\sqrt{3}&0&0&-1&0
\end{pmatrix}},\q
E_2={\tiny\begin{pmatrix}
0&0&\sqrt{3}&0&0\\
0&0&0&0&1\\
-\sqrt{3}&0&0&1&0\\
0&0&-1&0&0\\
0&-1&0&0&0
\end{pmatrix}},\q
E_3={\tiny\begin{pmatrix}
0&0&0&0&0\\
0&0&0&2&0\\
0&0&0&0&1\\
0&-2&0&0&0\\
0&0&-1&0&0
\end{pmatrix}},
\edm
or 
\begin{equation}
\label{eq:E-forms}
E_1={\sqrt{3}}e^{15}+e^{23}+e^{45},\q 
E_2={\sqrt{3}}e^{13}+e^{25}+e^{34},\q
E_3=2e^{24}+e^{35}
\end{equation}
when seen as elements of $\Lambda^2\R^5$.
\bigbreak
We shall say that \emph{a \NI\ structure ${\mathbb T}$,
  the corresponding torsion $T$ or the manifold carrying ${\mathbb T}$
has/is of type $\mathcal W$}\ if\ $\hodge T$ belongs to the
  irreducible module
$\mathcal W\set\Lambda^2$. The type can thus be
  $\Lambda^2_3\oplus\Lambda^2_7$ generically,\ $\Lambda^2_3$,
  $\Lambda^2_7$ which we shall both refer to as `pure' type, or
  $\{0\}$ i.e.~the torsion-free case.

\section{The  characteristic connection of a $5$-dimensional Lie group}
\noindent
We begin by taking a generic Lie group of dimension five and write its
structure equations
\begin{equation}
         \label{eq:generalLA}
          \left\{
            \begin{matrix}
              d\1  & = & b_1e^{12}    & + &  \ldots & + & b_{10}e^{45}\\
              d\2  & = & b_{11}e^{12} & + &  \ldots & + & b_{20}e^{45}\\
              d\3  & = & \hdotsfor[5]{5}\\
	 d\4  & = & \hdotsfor[5]{5}\\
              d\5  & = & b_{41}e^{12} & + &  \ldots & + & b_{50}e^{45}.
            \end{matrix}
          \right.
\end{equation}
Taking $(e_1,\ldots e_5)$ to be the adapted 
frame means fixing the $\SO(3)$ structure and varying the Lie algebra 
equations in terms of real numbers $b_\alpha,\ \alpha=1,\ldots, 50.$

Polarising the expression \eqref{eq:t_ijk} gives the components of
${\mathbb T}$
\bdm
       \ba{c}
       t_{111}=-1,\q t_{122}=1,\q t_{144}=1,\q t_{133}=-\tfrac 12,\q
       t_{155}=-\tfrac 12,\\[3pt]
       t_{433}=- \tfrac {\sqrt{3}}2,\q t_{455}=\tfrac {\sqrt{3}}2,\q
       t_{235}=\tfrac {\sqrt{3}}2.
       \ea
\edm
\begin{lemma}
The $\SO(3)$-tensor ${\mathbb T}$ is \NI\ if and
only if the structure equations satisfy the set of linear relations
\eqref{eq:conditions} below.
\end{lemma}
\begin{proof}
Allow $X$ to be a generic linear combination of the basis $(e_i)$ of the Lie
algebra $\l$ of $L$, i.e.~$X=\sum_1^5 \lambda_j e_j$. 
Imposing \eqref{eq:ni} means solving the equation
$(\nabla_X{\mathbb T})(X,X,X) =0$ for all
values of the $\lambda_j$'s.

Consider euristically $X=\lambda \e3 +\mu \e4$ and
let $Y=\na_XX$. Then
\eqref{eq:ni} becomes
\bdm\ba{rcl}
0&\!\!=\!\!&\!\! \tfrac 13(\nabla_X{\mathbb T})(X,X,X)={\mathbb T}(Y,X,X)=
\lambda^2{\mathbb T}(Y,\e3,\e3)+2\lambda\mu\,{\mathbb T}(Y,\e3,\e4)+
\mu^2{\mathbb T}(Y,\e4,\e4)\\[3pt]
      &\!\!=\!\!&\!\! \lambda^2(-\tfrac 12 g(Y,\e1)-\tfrac{\sqrt{3}}2
g(Y,\e4))+2\lambda\mu\,\tfrac {\sqrt{3}}2 g(Y,\e3)+\mu^2 g(Y,\e1)\\[3pt]
      &\!\!=\!\!&\!\! 
\lambda^2\bigl(\tfrac 12 g([X,\e1],X)+\tfrac{\sqrt{3}}2
g([X,\e4],X)\bigr) 
-{\sqrt{3}}\lambda\mu\,g([X,\e3],X)
-\mu^2 g([X,\e1],X),
\ea
\edm
eventually resulting in a fourth order polynomial $\sum_{i+j=4}
\lambda^i\mu^jP_{ij}$. Given that $\lambda,\mu$ are
arbitrary, the conditions translate into $P_{ij}=0$ for all
$i,j=0,\ldots, 4$. Explicitly
\bdm
\ba{c}
P_{40}=\tfrac 12 (b_{22}-{\sqrt{3}}b_{28})=0,\q P_{04}=-b_{33}=0,\q
P_{22}=\tfrac 12 b_{33}-b_{22}+{\sqrt{3}}b_{28}=0\\[3pt]
P_{31}=\tfrac 12(b_{32}+b_{23}-{\sqrt{3}}b_{38})=0, \q
P_{13}=-(b_{23}+b_{32})+{\sqrt{3}}b_{38}=0.
\ea
\edm
The story is similar when all $\lambda_i$'s are
present. Passing from linear combinations $X$ of $p$
elements to $p+1$ does not increase dramatically the complexity, since
most of the new information is trivial by previous
relations, justifying the fact that the system is linear.
The set of conditions imposed by \eqref{eq:ni} reads thus:
\begin{equation}
\label{eq:conditions}
\allowdisplaybreaks
\ba{c}
b_1=b_{11}=b_3=b_{33}=0,\ b_{20}= - b_{37},\ b_{13}+b_{31}=0,\\[3pt]
b_2={\sqrt{3}}(b_{23}+b_{8}),\ b_4={\sqrt{3}}(-b_{43}+b_{10}),\
b_{22}= {\sqrt{3}}b_{28},\
b_{44}={\sqrt{3}}b_{50},\\[3pt]
b_{21}+b_{12}={\sqrt{3}}b_{17},\ b_{14}+b_{41}=
{\sqrt{3}}b_{15},\ b_4={\sqrt{3}}(b_5-b_{21}),\ b_2={\sqrt{3}}(b_7-b_{41}),\\[3pt]
2b_{22}+{\sqrt{3}}b_{16}=2{\sqrt{3}}(b_{19}+b_{27}),\ 2b_{44}-{\sqrt{3}}b_{16}=
2{\sqrt{3}}(b_{45}-b_{19}),\\[3pt]
2b_{29}+b_{17}=b_{26}+b_{18},\ 2b_{29}+b_{40}=b_{26}+b_{35},\ 2b_{49}-b_{15} - b_{37}=b_{46},\\[3pt]
b_{28}+b_{50}=b_{45}+b_{27},\ b_{24}+b_{42}={\sqrt{3}}(b_{25}+b_{47}),\
       b_{48}-b_{30}=b_{47}-b_{25},\\[3pt]
2(b_{24}+b_9)=2{\sqrt{3}}b_{25}+b_{13}+b_6,\
2(b_{42}-b_9)=2{\sqrt{3}}b_{47}-(b_{13}+b_6),\\[3pt]
2(b_{39}+b_{30}-b_{25})=b_{36},\  b_{35}+b_{17}=b_{40}+b_{18},\ 
2(b_{48}+b_{39}-b_{47})=b_{36},\\[3pt]
b_{38}= -b_{15},\ {\sqrt{3}}(b_{40}+b_{18}-b_{35})=b_{21}+b_{12},\ 
b_{32} = - b_{23} - \sqrt{3} b_{15},\\[3pt]
\ea
\end{equation}
easily handled by computer programs.
\end{proof}
\bigbreak
Further constraints on the coefficients in \eqref{eq:generalLA} 
derive from $d^2=0$, but due to computational complexity we
reserve the Jacobi identity to when strictly necessary, typically at
the very end of the classifying process.
\bigbreak
Since $\tilde \na$ preserves the metric and
the tensor ${\mathbb T}$, it is clear that
\begin{center}
$\bg
g(\tilde \na_XY,Z)+g(\tilde \na_XZ,Y)=0,\\
{\mathbb T}(\tilde\na_XY,Z,W)+{\mathbb T}(\tilde \na_XZ,Y,W)+
{\mathbb T}(\tilde \na_XW,Y,Z)=0.
\eg $
\end{center}
The characteristic torsion \eqref{eq:torsion} of the structure is then
      given by
\bdm
\ba{rl}
T=&
(b_{43}-b_{10}+b_{12})e^{123}-b_6e^{124}+({\sqrt{3}}b_{15}-b_7)e^{125}+\\[4pt]
&
      ({\sqrt{3}}b_{15}-b_8)e^{134}+(b_{24}-{\sqrt{3}}b_{47}-
\tfrac 12 b_{13}-\tfrac 12 b_6)e^{135}+\\[4pt]
&
({\sqrt{3}}b_{40}-b_{10})e^{145}+(2b_{29}-b_{17}-b_{35})e^{234}+
(b_{28}-b_{50}-b_{19})e^{235}+\\[4pt]
& (b_{37}-b_{15}-2b_{49})e^{245}+(\tfrac {\sqrt{3}}6 b_{13}+\tfrac {\sqrt{3}}6
b_6-\tfrac {\sqrt{3}}3 b_9-\tfrac {\sqrt{3}}3 b_{24}+b_{39}-b_{47})e^{345}.
\ea
\edm

\noindent
The $SO(3)$-connection $\tilde\nabla$ can be thus reconstructed, and is
written here as an
$\so(3)$-valued form
$$
\Gamma = \gamma^1 E_1 + \gamma^2 E_2 + \gamma^3 E_3
$$ where
\begin{equation}
\label{eq:gamma-connection}
\gamma^1 = \Gamma^2_3, \qq \gamma^2 = \Gamma^2_5, \qq \gamma^3 = \Gamma^3_5
\end{equation}
and
\bdm
\ba{l}
    \Gamma^2_5 = \tfrac 1{\sqrt{3}} \Gamma^1_3 = \Gamma^3_4 = -( b_{23}
+  b_8) e^1 -
    b_{17} e^2 -  b_{28}e^3 + b_{15} e^4 -
    b_{47} e^5,\\[5pt]
\Gamma^2_3 = \tfrac 1{\sqrt{3}} \Gamma^1_5 =   \Gamma^4_5 = ( b_{43}
-  b_{10}) e^1 -
b_{15} e^2 + \tfrac 1{\sqrt{3}}(- b_9 +\tfrac 12 b_{13} + \tfrac 12
b_6 - b_{24}) e^3
-  b_{40} e^4 -  b_{50} e^5,\\[5pt]
\Gamma^3_5 =\tfrac 12\Gamma^2_4 =   -\tfrac 12(b_6 + b_{13}) e^1 +
( b_{45}-  b_{50} -  b_{19}) e^2 - b_{29} e^3 + ( b_{47}-  b_{39} -
b_{48} )e^4 -  b_{49} e^5
\ea
\edm
are the connection $1$-forms.
All $\Gamma^i_j$'s vanish precisely when
\begin{equation}
\label{eq:flatness}
\ba{c}
b_{43} = b_{10}, b_{23} = - b_8,  b_{15}= b_{17} = b_{29}  =
b_{40} = b_{49} =0,\\[5pt]
b_{28} = b_{47} = b_{50} =0, b_{45} = b_{19}, b_{48} = - b_{39},
b_{13} = - b_6,  b_{24} = - b_9.
\ea
\end{equation}
\bigbreak
Now the curvature $K=r^1E_1+r^2E_2+r^3E_3\in \Lambda^2\l^*\otimes
\so(3)$ is determined by
\bdm
\ba{c}
r^1=d\gamma_1+\gamma_2\w\gamma_3,\q 
r^2=d\gamma_2+\gamma_3\w\gamma_1,\q 
r^3=d\gamma_3+\gamma_1\w\gamma_2,
\ea 
\edm
and $r^1=K^2_3, r^2=K^2_5, r^3=K^3_5$ are the non-zero components of $K$.
If $T\equiv 0$ then $K$  is
determined by a constant $F\in\R$, since \cite{Bobienski-N:SO(3)}
\bdm
r^j=FE_j\qq j=1,2,3
\edm
and the $\SO(3)$ structure is locally isometric to that of a symmetric space
$Q/\SO(3)$, with
\begin{equation}
\label{eq:F}
Q=\left\{
\ba{ll}
\SO(3)\times_\rho\R^5 & \q\textrm{if}\ F=0\\
\SL(3,\R)          & \q\textrm{if}\ F<0\\
\SU(3)             & \q\textrm{if}\ F>0.
\ea
\right.
\end{equation}
\bigbreak
A purely formal but still suggestive
relation between near integrability and lower dimensional geometry
is the fact that
the number of parameters $b_\alpha$ upon which the
Riemannian connection of a \NI\ structure
depends, is precisely half of the total, exactly
as for self-dual connections on a four-manifold
\cite{Nurowski:distinguished-dimensions}. 
Whether this has to 
do with the universal covering $\SU(2)\to\SO(3)$ remains to be seen,
though the fact that the Berger sphere $\SO(5)/\SO(3)$ is
diffeomorphic to $\Lie{Sp}(2)/\SU(2)$ is a clear indication, and
will be dealt with elsewhere by the authors. But here
is an
\begin{ex*}
Consider the Lie algebra $\mathfrak l$ with structure equations 
\begin{equation}\label{ex:hypo}
\left\{
\ba{l}
de^1= b_7 e^{25} -b_{37} e^{35} -b_{45} e^{45}, \\
de^2= -b_7 e^{15} +b_{45} e^{35}  -b_{37} e^{45}, \\
de^3= b_{37} e^{15} -b_{45} e^{25}+\frac{(b_{37}^2+b_{45}^2)}{b_7} e^{45},\\
de^4= b_{45} e^{15} +b_{37} e^{25}  - \frac{(b_{37}^2+b_{45}^2)}{b_7}  e^{35},\\ 
de^5= b_7 e^{12}  -b_{37}e^{24} -b_{45} e^{14} -b_{37} e^{13} +
b_{45} e^{23} + \frac{(b_{37}^2+b_{45}^2)}{b_7}  e^{34},
\ea
\right.
\end{equation}
and the \NI\ structure  whose associated three-form $T$
is given by  
$$
T = -b_7 e^{125}  +b_{37} e^{135} +b_{45} e^{145}- b_{45} e^{235} +
b_{37} e^{245} - \tfrac
{(b_{37}^2+b_{45}^2)}{b_7} e^{345}.
$$
First of all this torsion is harmonic. But more crucially, the
associated Lie group $L$ is thus endowed 
with an $\SU(2)$ structure defined by 
 $$
 \alpha = e^5,\q \omega_1=e^{12}+e^{34},\ \omega_2=e^{13}+e^{42},\ 
\omega_3=e^{14}+e^{23}.
 $$
To be precise, this is actually hypo 
\cite{Conti-S:KS-in-dim5} as 
\[
d\omega_1 = 0,\q d (\omega_2 \wedge \alpha) = d (\omega_3 \wedge
\alpha) =0.
\]
\noindent
Examples of this sort are noteworthy because, being hypo, they induce
local integrable $\SU(3)$ structures on $L\times
\R$. 
Moreover, they yield half-flat geometries, which also evolve in one
 dimension higher, but to holonomy
$\G$. In fact, \eqref{ex:hypo} is part of a
more general family, as shown by
\begin{proposition}
Let $L$ be a Lie group whose Maurer-Cartan structure 
\eqref{eq:generalLA} satisfies \eqref{eq:flatness} plus 
$ 
b_6 = b_8 = b_{12} = b_{35} =0, b_{10} = - b_{45},  b_9 = - b_{37}, b_{39} = (b_{37}^2 +
b_{45}^2)/b_7.
$
Then for some closed $1$-form $e^6$, the invariant $\SU(3)$ structure on 
$L\times \R$
$$
\omega = \omega_3 + \alpha \wedge e^6,\qq
 \psi_+ = \omega_2 \wedge \alpha - \omega_1 \wedge e^6 
$$
is half-flat.
\end{proposition}
\begin{proof}
Almost immediate,  once one recalls that half-flat is a fancy name for
the closure of $\psi_+$ and $\omega\w\omega$.
\end{proof}
\end{ex*}

\section{Classification}

\noindent
The action of the endomorphisms ${\mathbb T}_X={\mathbb T}(X,\cdot,\cdot)$
indicates that whilst ${\mathbb T}_{\e1},
{\mathbb T}_{\e2}$ and ${\mathbb T}_{\e4}$ preserve the decomposition
\begin{equation}
\label{eq:dec}
\l = \h \oplus \p=\text{span}\{e_1,e_2,e_4\}\oplus\ \text{span}\{e_3,e_5\},
\end{equation}
${\mathbb T}_{\e3}, {\mathbb T}_{\e5}$ on the contrary do not. This
reflects the
irreducibility of the rotational action.

We shall capture the Lie algebras $\l$ for which the
$SO(3)$-connection satisfies the conditions:
\begin{equation} \label{SO(3)conditions}
\begin{array}{l}
\tilde\nabla_X {\gh h} \subseteq {\gh h}, \quad \tilde\nabla_X {\gh p}
\subseteq {\gh p}, \quad \forall X \in {\gh h}\\
\tilde\nabla_Y {\gh h} \subseteq {\gh p},  \quad \tilde\nabla_Y {\gh p}
\subseteq {\gh h}, \quad \forall Y \in {\gh p}.
\end{array}
\end{equation}
Otherwise said, $\h$-derivatives preserve the splitting
\eqref{eq:dec}, but differentiation in the $\p$-direction swaps the
subspaces.  The geometric motivation of this 
condition are readily explained.

In order to classify $5$-dimensional Lie
groups $L$ with an invariant $SO(3)$ structure $(g, {\mathbb T})$ whose
    characteristic connection $\tilde\nabla$ satisfies
\eqref{SO(3)conditions}, it is quite useful to distinguish whether
$({\mathfrak l}, {\mathfrak h})$ is a  symmetric pair or not.
A Lie algebra ${\mathfrak l}= {\mathfrak m}\oplus \gh{q}$, or better 
$(\l,\gh{m})$ is called a \emph{symmetric pair}
when
$$
[ {\mathfrak m}, {\mathfrak m} ] \subseteq {\mathfrak m}, \quad
[{\mathfrak m}, {\mathfrak q}] \subseteq {\mathfrak q}, \quad
[{\mathfrak q}, {\mathfrak q} ] \subseteq {\mathfrak m}.
$$

\noindent
If $H$ is the simply-connected Lie subgroup of
$L$ with Lie algebra $\mathfrak h$, we have a Riemannian
fibration of Lie groups 
\begin{equation}
\label{eq:fibration}
H \longrightarrow L \longrightarrow L/H
\end{equation}
whose fibres $H$ are totally geodesic submanifolds of 
$L$ and whose base is a $2$-dimensional symmetric space.
Therefore $L/H$ is a finite quotient of $\R^2$, the Riemann sphere
or the Poincar\'e disc according to the curvature. 
Since $\p$ is indeed the tangent $2$-plane spanned by $e_3, e_5$, the
curvature of $L/H$ is given by
\begin{equation}
\label{eq:k(p)}
k(\p) = - g(e_5, \big[ [e_3, e_5]_{\mathfrak h}, e_3 \big]_{\mathfrak p})
 = b_9 ( \tfrac 12 b_6 +\tfrac 12 b_{13} -b_9 -{\sqrt{3}}b_{47}) -
     b_{19} b_{45} + b_{39} b_{48},
\end{equation}
\noindent 
and its  sign decides which of the space forms one is looking at.

\begin{lemma}
If ${\gh p}$ is $\ad({\gh h})$-invariant then ${\gh h}$ is a Lie
subalgebra of ${\gh l}$ and $[{\gh p}, {\gh p}]
\subseteq {\gh h}$.
\end{lemma}
\begin{proof}

If  $\p$ is $\ad(\h)$-invariant, the structure coefficients
\bdm
b_2, b_4, b_5,b_7,b_8,b_{10},\  b_{12}, b_{14},b_{15} ,b_{17},
b_{18},\ b_{32},b_{34},b_{35},b_{37},
b_{38},b_{40}
\edm
all vanish. As a consequence of \eqref{eq:conditions} also
$b_{23},b_{43},b_{21},b_{41}$ are zero.
     So the Maurer--Cartan system \eqref{eq:generalLA} simplifies to
\bdm
\left\{
\ba{ll}
d\1 =& b_6e^{24} + b_9e^{35}\\[3pt]
d\2 =& b_{13}e^{14}+2(b_{50}+b_{19}-b_{45})e^{24}+b_{19}e^{35}\\[3pt]
d\3 =& {\sqrt{3}}b_{28}e^{13}+b_{24}e^{15}+\tfrac
{\sqrt{3}}3(b_{9}+b_{24}-\tfrac
12 (b_{13}+b_{6}))e^{23}+2b_{29}e^{24}+
(b_{50}+b_{28}-b_{45})e^{25}+\\[3pt]
& b_{28}e^{34}+b_{29}e^{35}+\tfrac {\sqrt{3}}3(b_{9}+ b_{24}-\tfrac 12 (b_{13}+
b_{6})+{\sqrt{3}}b_{48}-{\sqrt{3}}b_{47})e^{45}\\[3pt]
d\4 =& -b_{13}e^{12}+2(b_{39}+b_{48}-b_{47})e^{24}+b_{39}e^{35}\\[3pt]
d\5 =&  (b_9+{\sqrt{3}}b_{47}-\tfrac 12 b_{13}-\tfrac 12
b_{6})e^{13}+{\sqrt{3}}b_{50}e^{15}+b_{45}e^{23}+2b_{49}e^{24}+
b_{47}e^{25}+b_{48}e^{34}+\\[3pt]
& b_{49}e^{35}+b_{50}e^{45}.
\ea
\right.
\edm
Using the vanishing  of the  coefficients  of $e^{234},
e^{245}$ in
$d(de^i), i = 1, 2, 4,$  and of $e^{124},
e^{135}, e^{235}, e^{345}$
for $i = 3, 5$, then either
\begin{enumerate}
\item[$\cdot$] $b_{29} = b_{49} = 0$ (implying $[{\gh p}, {\gh
p}] \subseteq {\gh h}$), or
\item[$\cdot$] $b_{29}^2 + b_{49}^2 \neq 0$ and
$
b_9 = b_{19} = b_{39} =0, b_{13} = - b_6, b_{45} = - b_{28}, b_{48} = b_{47},
b_{50} = - b_{28}
$
(which force ${\gh p}$ to become a Lie subalgebra of ${\gh l}$).
\end{enumerate}
But in the latter case, the Jacobi equation ends up
annihilating both $b_{49}$ and $b_{29}$,  
contradicting the assumption.
\end{proof}

\begin{corollary}
Demanding ${\gh p}$ to be $\ad({\gh h})$-invariant
      is equivalent to  the pair $({\mathfrak l}, {\mathfrak h})$
being symmetric.
\end{corollary}
\begin{proof}
In terms of the structure equations $\tilde\nabla$
satisfies \eqref{SO(3)conditions} if and only
if
$$
b_{43} = b_{10}, b_{23} = - b_8,  b_{15}= b_{17} = b_{29} = b_{37} =
b_{40} = b_{49} =0.
$$
If, in addition, $b_7 = b_8 = b_{10} = b_{12} = b_{35} =0$ then
$[\p,\h]\set\p$ and ${\gh h}$ is a
Lie subalgebra of ${\mathfrak l}$.
\end{proof}
\bigbreak
%
%
The algebraic structures appearing later will often be \emph{solvable} Lie
algebras $\g$, whose derived series
\bdm
\g\supseteq\g^1\supseteq\g^2\supseteq\ldots\supseteq\g^q=0
\edm
collapses at some stage $q\in\bb{N}$ called the step length.
Each ideal is defined by bracketing the
preceding one in the sequence with itself, $\g^{i+1}=[\g^i,\g^i]$,
beginning from $\g^0=\g$. The first subspace $\g^1=[\g,\g]$ is 
the commutator and such a $\g$ is referred to as $q$-step (solvable).
\begin{exs*}
There is only one solvable (non-nilpotent) Lie algebra of dimension
two, modulo isomorphisms:

\centerline{$\s_2 =(0,e^{12})$,}

\noindent 
and similarly, only one $3$-dimensional $2$-step one, namely

\centerline{$\s_3 =(0,0,e^{13})$.}

This notation for Lie algebras expresses the bracket of $\g$ via
the exterior differential of $\g^*$, so in the latter example one should understand the
basis $e_1,e_2,e_3$ of $\s_3$ to satisfy $de^1=de^2=0, de^3=e^{13}$,
or equivalently $[e_1,e_2]=[e_2,e_3]=0,[e_1,e_3]=-e_3$.
\end{exs*}

Solvable Lie algebras and groups  are intimately linked to non-positive curvature, as
predicted by Alekseevski\u{\i}
\cite{Alekseevsky:solvable}.
This phenomenon manifests itself patently in the coming sections, goal
of which is to prove the
\begin{theorem}
\label{thm:main}
If the characteristic connection behaves as in
\eqref{SO(3)conditions}, the torsion is always coclosed $d\hodge
T=0$. Moreover

%
%
%
%

\noindent
a) if $(\l, \h)$ is symmetric, $\l$ is either solvable or
isomorphic to
\bdm\ba{l}
\mathfrak{l_1}=\so(3)\oplus\R^2\\[3pt]
\mathfrak{l_2}= (2e^{12}, e^{14} ,  e^{15} -  e^{23},   2e^{24},  
e^{25} + e^{34}),\\[3pt]
\mathfrak{l_3}= \so(3)\oplus\s_2
\ea\edm
\noindent
b) If $(\l, \h)$ is not symmetric instead, then $\tilde \nabla \equiv 0$ 
and  $\l$ is isomorphic to $\so(3)\oplus\R^2$. 
\end{theorem}
\noindent
The appearance of algebras isomorphic to $\so(3)\oplus\R^2$ 
in distinct contexts depends
upon the underlying $\SO(3)$ structure, which changes in \emph{a)} and
\emph{b)}.
The proof of this statement is scattered over
sections \ref{sec:symmetric pair} and \ref{sec:non-symmetric}.
\bigbreak
The type of $\SO(3)$ geometry is determined by the 
presence of intrinsic torsion in the modules
\eqref{eq:torsion-types}. Because the two-forms $E_i$'s are eigenvectors of the
endomorphism $\hat {\mathbb T}$ of \eqref{eq:torsion-pairing} with a known
eigenvalue, one has that $\hodge
T\in\Lambda^2_7$ if and only if 
\begin{equation}
\label{eq:pure}
T\w E_i=0,\ i=1,2,3.
\end{equation}
The Lie algebras of type $\Lambda^2_7$ 
are therefore those for which this linear system 
in the $b_\alpha$'s holds.
They will be highlighted in the rest of the paper.
\bigbreak
To conclude the section we recall two basic equations
\cite{Ivanov-P:vanishing} 
for metric connections with torsion $\tilde\na$. One  
relates the covariant and exterior derivatives of $T$ 
\begin{equation}
\label{eq:nablaT-dT}
dT(x,y,z,w)=\gh{S}_{x,y,z}(\tilde\na_xT)(y,z,w)-(\tilde\na_wT)(x,y,z)
+2\sigma_T(x,y,z,w),
\end{equation}
where $\gh{S}$ is the cyclic sum.  
The four-form 
\begin{equation}
\label{eq:sigmaT}
2\sigma_T=\tsum_1^5 (e_i\hook T)^2
\end{equation}
also shows up 
in a second equation, that is the Bianchi identity for the curvature
$\tilde R$ of $\tilde\na$
\[
\gh{S}_{x,y,z}\tilde R(x,y,z,w)=dT(x,y,z,w)-\sigma_T(x,y,z,w)
+(\tilde\na_wT)(x,y,z).
\]
\begin{corollary}
\label{cor:IP}
\cite{Ivanov-P:vanishing} 
When the torsion is $\tilde\na$-parallel then $dT=2\sigma_T$.

If additionally the manifold  happens to be flat $\tilde R=0$, then it has closed
torsion.\qed
\end{corollary}
For the significance of $\sigma_T$ see
 for instance \cite{Agricola-F:holonomy}.

\section{$({\mathfrak l}, {\mathfrak h})$ symmetric pair}
\label{sec:symmetric pair}

\noindent
The case when $({\mathfrak l}, {\mathfrak h})$ is a symmetric pair
corresponds to the Levi-Civita connection
satisfying a bunch of relations similar to  \eqref{SO(3)conditions}:

\begin{proposition}
\label{prop:LC conditions}
The Levi-Civita connection $\nabla$ satisfies
\bdm 
\begin{array}{l}
\nabla_X {\gh h} \subseteq {\gh h}, \quad \nabla_X
{\gh p} \subseteq {\gh p}, \quad \forall X \in {\gh h}\\
\nabla_Y {\gh h} \subseteq {\gh p},  \quad \nabla_X
{\gh p} \subseteq {\gh h}, \quad \forall Y \in {\gh p}
\end{array}
\edm
if and only if ${\gh p}$ is $\ad({\gh h})$-invariant and ${\gh
h}$ is a Lie subalgebra of ${\mathfrak l}$.
\end{proposition}

\begin{proof} The conditions are equivalent to
the vanishing of the coefficients $b_7$, $ b_8$, $ b_{10}$, $
b_{12}$, $ b_{15}$, $ b_{17}$, $ b_{23}$, $ b_{29}$, $ b_{35}$,
$b_{37}$, $ b_{40}$, $b_{43}$, $ b_{49}$ in \eqref{eq:generalLA}.
\end{proof}

\noindent
The first relation is an old acquaintance, for
$0=(\na_\h\h)_\p$ is the second
fundamental form of $H\sset L$.
\bigbreak
Imposing \eqref{SO(3)conditions} reduces the characteristic connection
to the treatable
\begin{equation}
\label{eq:char-connection}
\ba{l}
 \tilde\nabla e_1 =- {\sqrt{3}}(b_{28} e_3 + A e_5)\otimes 
e^3 +{\sqrt{3}}(b_{47} e_3 + b_{50} e_5)\otimes e^5\\[4pt]
 \tilde\nabla e_2 = 2 B e_4\otimes e^1+2  C e_4\otimes e^2 
+(A e_3 + b_{28} e_5)\otimes e^3+2 D e_4\otimes e^4+\\
  \ph{ \tilde\nabla e_1 =-}\q (b_{50}e_3 + b_{47} e_5)\otimes
  e^5,\\[4pt]
 \tilde\nabla e_3 = B e_5\otimes e^1+C e_5\otimes 
e^2+\bigl(b_{28}(e_4-{\sqrt{3}}e_1) - A e_2\bigr)\otimes e^3+D 
e_5\otimes e^4+\\
  \ph{ \tilde\nabla e_1 =-}\q 
\bigl(b_{47}(e_4-{\sqrt{3}}e_1)- b_{50} e_2\bigr)\otimes e^5,\\[4pt]
 \tilde\nabla e_4 = - 2 B e_2\otimes e^1- 2 C e_2\otimes 
e^2+(-b_{28} e_3 + A e_5)\otimes e^3- 2 D e_2\otimes e^4+\\
  \ph{ \tilde\nabla e_1 =-}\q (b_{50} e_5 - b_{47} e_3)\otimes e^5,\\[4pt]
 \tilde\nabla e_5 =  - B e_3\otimes e^1- C e_3\otimes e^2-\bigl( 
b_{28} e_2 + A(e_4-{\sqrt{3}}e_1)\bigr)\otimes e^3- D e_3\otimes 
e^4-\\
  \ph{ \tilde\nabla e_1 =-}\q 
\bigl(b_{47} e_2 +b_{50}(e_4+{\sqrt{3}}e_1)\bigr)\otimes e^5,\\
\ea
\end{equation}
where capitals are merely used to abbreviate the constants
\bdm
\ba{c}
 B=\tfrac
12(b_6+b_{13}),\q A=\tfrac{{\sqrt{3}}}{3}(b_9+b_{24}-B),\q 
C=b_{50} - b_{45} + b_{19},\q  D=b_{39} + b_{48} - b_{47}.
\ea
\edm

\begin{corollary}
The splitting $\h\oplus\p$ is curvature-invariant:
$$
R_{XY}(\h)\set\h,\ R_{XY}(\p)\set\p\qq \forall X,Y\in \l.
$$
\end{corollary}
\begin{proof}
   From \eqref{SO(3)conditions} one induces easily that
$$
\begin{array}{l}
R_{{\mathfrak h} {\mathfrak h}} {\mathfrak h},\ R_{{\mathfrak h}
{\mathfrak p}} {\mathfrak h},\ R_{{\mathfrak p} {\mathfrak p}}
{\mathfrak h} \in {\mathfrak  h},\\
R_{{\mathfrak h} {\mathfrak p}} {\mathfrak p},\ R_{{\mathfrak h}
{\mathfrak h}} {\mathfrak p},\ R_{{\mathfrak p} {\mathfrak p}}
{\mathfrak p} \in {\mathfrak p}.
\end{array}
$$
\noindent
Using the symmetries of $R$, this means $R_{{\mathfrak h}  {\mathfrak h}  {\mathfrak h}
{\mathfrak p}}= R_{{\mathfrak h} {\mathfrak p} {\mathfrak p}
{\mathfrak p}} =0$, i.e.~whenever any three indexes in the 
components $R_{ijkl}$ denote vectors in $\p$ (or $\h$), the curvature is
zero. In other words the curvature $2$-form preserves \eqref{eq:dec}.
\end{proof}
\noindent
As a consequence, the $\SO(3)$-curvature $\tilde R$ fulfills
similar conditions 
$ \tilde R_{\h\h\h\p}=\tilde R_{{\mathfrak h}  {\mathfrak p} 
{\mathfrak h}  {\mathfrak h}}
= \tilde R_{{\mathfrak h} {\mathfrak p} {\mathfrak p}
{\mathfrak p}}=\tilde R_{\p\p\h\p}=0 $, despite it no
longer being a symmetric endomorphism of $\Lambda^2\l$.
\bigbreak
The vector space decomposition \eqref{eq:dec} transfers
also to the algebraic level as follows
\begin{corollary}
When $\p$ is Abelian ${\mathfrak l}$ is a semidirect sum
$${\mathfrak l} = {\mathfrak h} \oplus_{\alpha} {\mathfrak
     p},$$
where $\alpha:{\mathfrak h} \to {\mathfrak {Der}} ({\mathfrak
p})$ is a homomorphism of Lie algebras.
\end{corollary}
\begin{proof}
This descends  from the $\ad({\mathfrak h})$-invariancy of $\p$.
\end{proof}
\bigbreak
For better handling the discussion now divides into the mutually
exclusive cases
\begin{center}
$T=0$, \qq $dT=0\ (T\neq 0$, called `strong'), \q and \q $dT\neq 0$,
\end{center}
the simplest being the torsion-free one.


\begin{proposition} Let $(\l, \h)$ be  a symmetric pair with $T\equiv
  0$. Then
\begin{itemize}
\item[(i)] \cite{Bobienski-N:SO(3)} the $\SO(3)$ structure of $L$ is 
locally isometric to that of a
 symmetric space of non-positive curvature, hence either 
$\R^5$ or $\SL(3,\R)/\SO(3)$;
\item[(ii)] $\l$ is solvable;
\item[(iii)] the symmetric space $L/H$ is flat.
\end{itemize}
%
%
%
\end{proposition}
\begin{proof}

     The structure equations are
\bdm
\left\{
\ba{l}
d\1  =  b_9e^{35}\\[3pt]
d\2  =
b_{13}e^{14}+2(b_{50}+b_{19}-b_{45})e^{24}+(b_{28}-b_{50})e^{35}\\[3pt]
d\3  =  {\sqrt{3}}b_{28}e^{13}+(\tfrac 12
b_{13}+{\sqrt{3}}b_{47})e^{15}+(-b_{47}+\tfrac {\sqrt{3}}3
b_{9})e^{23}+	 (b_{50}+b_{28}-b_{45})e^{25}+\\[3pt]
\phantom{d\3  =\  }b_{28}e^{34}+(\tfrac {\sqrt{3}}3
b_{9}+b_{48})e^{45}\\[3pt]
d\4  = -b_{13}e^{12}+2(b_{39}+b_{48}-b_{47})e^{24}+b_{39}e^{35}\\[3pt]
d\5  =
(b_9+2{\sqrt{3}}b_{47}-b_{24})e^{13}+{\sqrt{3}}b_{50}e^{15}+
b_{45}e^{23}+b_{47}e^{25}+b_{48}e^{34}+b_{50}e^{45}.
\ea
\right.
\edm
The inspection of the Jacobi identity tells that there are four
solutions for which $\na=\tilde\na$.
Only the non-zero coefficients  are indicated
and serve to distinguish the Lie algebras:
\begin{enumerate}

\item $b_{28}=b_{50} =a \neq 0$\\[3pt]
$(0, 2 a  e^{24}, a ({\sqrt{3}} e^{13} +   e^{34}
+ 2 e^{25}), 0, a ({\sqrt{3}} e^{15} + e^{45}))$
\smallbreak

\item $b_{28} = b_{50} =\tfrac 12 (b_{45}^2+b_{48}^2)/b_{45}, b_{45} =a
\neq 0, b_{48} =b$\\[3pt]
$(0,   \frac{- a^2 +b^2}{a} e^{24},
      b  e^{45} +  \frac{a^2 +b^2}{2a} (e^{34} +
{\sqrt{3}}  e^{13})
+  \frac{b ^2}{a} e^{25},  2 b  e^{24},
\frac{ a^2 +b^2}{2a} (e^{45} +{\sqrt{3}} e^{15}) + b e^{34}
+ a  e^{23})$.
\smallbreak

\item $b_{50} = -\tfrac 12 b_{19}, b_{28} = \tfrac 12 b_{19},  b_{39} =
2b_{47}, b_{24} =  {\sqrt{3}} b_{47}, b_{19} = a, b_{47} =b, a^2 +
b^2 \neq 0$\\[3pt]
$(0, a  (e^{24} +  e^{35}),
\tfrac 12 a ( {\sqrt{3}}   e^{13} + e^{34}) + b ({\sqrt{3}} e^{15} + e^{23}),
b (2 e^{24}  + 2  e^{35}),
b ( {\sqrt{3}} e^{13} + e^{25})  - \tfrac 12 a (
{\sqrt{3}} e^{15} + e^{45}))$.
\smallbreak

\item $b_{13} = 2b_{24} =  a \neq 0$\\[3pt]
$(0, a e^{14},  \tfrac 12 a e^{15}, -a e^{12},
      -\tfrac 12 a e^{13})$.

\end{enumerate}

\noindent
The first three instances are
$3$-step solvable Lie algebras with $3$-dimensional commutator
${\mathfrak l}^1$ and  ${\mathfrak h} \cong \R \oplus {\mathfrak
s}_2$.  The Lie algebra (1) is isomorphic to
(2) with $b =0$. The Lie algebra $(4)$ is
$2$-step solvable with
$4$-dimensional commutator ${\mathfrak l}^1$ and $\h$ Abelian.
Instead $\mathfrak p$ is an
Abelian subalgebra only for the Lie algebras (1), (2) and (4), so 
the curvature formula \eqref{eq:k(p)} implies
 sectional flatness. In
the remaining case (where $\mathfrak p$ is not Abelian) we still
have $k(\p) =0$.
The constant $F$ of \eqref{eq:F} is negative for (1)-(3) and
vanishes for (4), so \cite{Bobienski-N:SO(3)} implies that the first three are
isometric to $\SL(3,\C)/\SO(3)$, the last to $\R^5$.
\end{proof}

The Lie algebras listed (in roman numerals) in the sequel have
simpler structure equations than those generated in the proofs. They are
attained by standard changes of bases
-- omitted not to bore the reader senseless -- and may not be optimal: 
for example $\sqrt{3}e^1+e^4$ could reasonably be a basic $1$-form, as variously
hinted by, e.g., \eqref{eq:char-connection}. This depends
essentially on the choice of representation \eqref{eq:rep}.

Nevertheless, it should be all but clear that doing this will 
alter the $\SO(3)$ structure rather dramatically, so retaining in the 
proofs the original expressions, prior to any algebra isomorphism, allows to 
read off the geometry of interest.
In order to truly distinguish algebras up to $\SO(3)$
equivalence, Cartan-K\"ahler theory seems the only reasonable way. 
This was pursued in \cite{Bobienski-N:SO(3)}.
\bigbreak
%
%
\begin{lemma}
\label{lemma:symm-closed}
If $(\l, \h)$ is a symmetric pair with strong torsion, then
$\l$ must be isomorphic to one of
\begin{enumerate}

\item[I.] \q
$(0,  -2  e^{24}, (e^{13}+ e^{34}) + m ( e^{15}   + e^{23}),
-2 m e^{24}, m ( e^{13} + e^{25}) -  (e^{15} +e^{45}))$
\smallbreak

\item[II.] \q
$(e^{24}, e^{14}, 0, e^{12}, 0)$
\smallbreak

\item[III.] \q
$ ( e^{35}, 0,e^{45} +e^{23},  e^{24}, -e^{25} )$
\smallbreak

   \item[IV.] \q
   $( -c e^{35}, a e^{35}, c e^{15} - a  e^{25}  + b e^{45},   -b
   e^{35},   -c e^{13} + a e^{23}+ b  e^{34})$,
   \end{enumerate}
where $ m, a, b, c \in \R, a^2 +b^2 + c^2 \neq 0$.
The first two algebras give rise to a flat quotient $L/H$, whereas 
$k(\p)$ is negative for
III  and IV. 
The torsion is of pure type
$\Lambda^2_7$ in the family of Lie algebras IV when $c=0$, and always in I.
\end{lemma}

\begin{proof}

Suppose that $T$ be closed, and not identically zero.
After imposing the Jacobi identity one reduces to four instances:

\begin{enumerate}

\item $b_{50} = -b_{28} = - a \neq 0,   b_{24} =  {\sqrt{3}} b_{47} =
{\sqrt{3}} b$   \\[3pt]
$\vs(0,  -2 a e^{24},
       a ({\sqrt{3}} e^{13}+ e^{34}) + b ({\sqrt{3}}  e^{15}   +
e^{23}),    -2 b e^{24},
b ({\sqrt{3}}  e^{13} + e^{25}) - a ({\sqrt{3}}  e^{15} +
     e^{45}))$
with $T = 2(a e^{2}+be^{4})e^{35}$ and
$k(\p)=0$.
Example 6.3.3 in \cite{Bobienski-N:SO(3)} is
contained in this family of algebras by requiring $b\neq 0$.
\smallbreak

\item $ b_{6} = -b_{13} = a \neq 0$ \\[3pt]
$(- a e^{24},  a  e^{14}, 0,
      -a e^{12}, 0)$ with
$T = a e^{124}$ and
$k(\p)=0$.
\smallbreak

\item $ b_{47} = -\tfrac 13 {\sqrt{3}} b_{9} = -  a \neq 0$\\[3pt]
$ ( {\sqrt{3}} a e^{35}, 0,   2a e^{45}  +  a  e^{23},
       2a e^{24},   - a e^{25})$
with $T =   {\sqrt{3}} a e^{135} $ and
$k(\p)=-3 a^2$.
\smallbreak

\item $ b_{45} = b_{19}= a,  b_{39} = -b_{48}= - b,   b_{9} = -b_{24} = - c,
a^2 +b^2 + c^2 \neq 0$\\[3pt]
$( -c e^{35}, a e^{35}, c e^{15} - a  e^{25}  + b e^{45},   -b
e^{35},   -c e^{13} + a e^{23}+ b  e^{34})$ with\\[3pt]
$T =  (c e^{1}-a e^{2} +b e^{4})e^{35} $ and
$ k(\p)=  - a^2 - b^2 - c^2$.
\smallbreak

\end{enumerate}
The following observations guarantee that three cases are
algebraically distinct.
The first Lie algebra is $2$-step solvable with
$\dim\l^1=3$, has an Abelian ${\mathfrak p}$ and ${\mathfrak h} \cong \R
\oplus {\mathfrak s}_2$. The third one is $3$-step solvable with 
$4$-dimensional commutator.
Numbers (2), (4) are both essentially $\so(3)\oplus\R^2$, 
and account for $\gh{g_1}$ in theorem \ref{thm:main}.
\end{proof}
\bigbreak
%
%
\begin{lemma}
\label{lemma:symm-non-closed}
If $dT\neq 0$ and $(\l, \h)$ is a symmetric pair  then
$\l$ is one of
\begin{enumerate}

\item[I.] \q
$(  e^{24},  0,  - e^{23} ,  e^{24}, e^{25} +   e^{34})$
\smallbreak

\item[II.] \q
$( 6  e^{24}, 2 e^{14} ,  e^{15} -   {\sqrt{3}} e^{23} -   {\sqrt{3}}
e^{45},   -2e^{12},  -  e^{13} +  {\sqrt{3}} e^{25} +   {\sqrt{3}}
     e^{34})$
\smallbreak

\item[III.] \q
$\bigl( (2c+ {\sqrt{3}}(b^2+1) ) e^{24} + c e^{35},
      -2 b e^{24},
       - \tfrac 12  (b^2+1) e^{23}  - b  e^{25} -
b^2 e^{45},(1-b^2) e^{24},
      b e^{23} + \tfrac 12  (b^2+1) e^{25} +
e^{34}\bigr)$
\smallbreak

\item[IV.] \q
$\bigl(-2 (b^2+c^2) e^{24} - (b^2+c^2 +
1) e^{35},
       -2 b e^{24},  e^{15} - b e^{25} + c
e^{45},   2 c e^{24},   - e^{13} + b e^{23}  +
ce^{34}\bigr)$
\end{enumerate}
up to isomorphisms, where $  b, c \in \R$.
The surface $L/H$ has zero curvature in cases I, II and negative in
the others.
\end{lemma}
\begin{proof}

      If $d T \neq 0$ the following possibilities come out of $d^2=0$:

\begin{enumerate}

\item $b_{6} = 2 {\sqrt{3}} b_{47}=2{\sqrt{3}}a \neq 0$; \\[3pt]
$( 2 {\sqrt{3}} a e^{24},  0,
      -a e^{23}  - 2 ae^{45},
      -2 ae^{24},  a  e^{25})$  and
$$ T =  -2 {\sqrt{3}} a (e^{124} + e^{135}), \quad
dT= 12 a^2 e^{2345}.$$
\smallbreak

\item $  b_{48} = 2b_{47}, b_{6} = 2 {\sqrt{3}}b_{47}, b_{47}=a \neq 0$\\[3pt]
$(2 {\sqrt{3}} a e^{24},  0,
       -a e^{23} ,
2 a e^{24}, a e^{25} + 2 a e^{34})$ and
$$T =  -2{\sqrt{3}}a(e^{124} + e^{135}), \quad
dT= 12a^2 e^{2345}.$$
\smallbreak

\item $ b_{24} = \tfrac 12 b_{13}, b_{48} = \tfrac {\sqrt{3}}2 b_{13}, b_{6}
= 3b_{13}, b_{47} = \tfrac {\sqrt{3}}2 b_{13}, b_{13}  = a \neq 0$\\[3pt]
$( 3 a e^{24},        a  e^{14} ,
\tfrac 12 a e^{15} - \tfrac {\sqrt{3}}2 a e^{23} - \tfrac {\sqrt{3}}2 a
e^{45},   -a e^{12},
      - \tfrac 12  a e^{13} + \tfrac {\sqrt{3}}2 a e^{25} + \tfrac {\sqrt{3}}2
a e^{34})$ and
$$T =  -3a (e^{124} + e^{135}), \quad
dT = 9a^2 e^{2345}.$$
\smallbreak

\item $  b_9 = a, b_{45} = b, b_{48} = c\neq 0,  b_{6} =
\tfrac{2 a c+{\sqrt{3}} b^2+{\sqrt{3}} c^2}{c}, b_{47} =
\tfrac {b^2+c^2}{2c}   $\\[3pt]
$( \frac{2a c + {\sqrt{3}} b^2+ {\sqrt{3}} c^2}{c} e^{24} +
a e^{35},  -2 b e^{24},
       - \frac{b^2+c^2}{2c} e^{23}  - b  e^{25} -
\frac{b^2}{c}  e^{45},\frac{c^2-b^2}{c} e^{24},
      b e^{23} +  \frac{b^2+c^2}{2c}  e^{25} + c
e^{34})$ with
$$
T =  -\tfrac{2ac+ {\sqrt{3}} b^2+ {\sqrt{3}} c^2}{c} e^{124}
-\tfrac{a c+ {\sqrt{3}} b^2+ {\sqrt{3}} c^2}{c} e^{135},\q
dT = \tfrac{(2 a c+ {\sqrt{3}} b^2+ {\sqrt{3}} c^2)^2}{c^2} e^{2345}
$$
and $k(\p)=-a^2$.
The algebra ${\mathfrak l}$ is $3$-step solvable. For $b \neq 0$ the
dimension of ${\mathfrak l}^1 $ is four. Further asking $a \neq 0$
recovers \cite[Example  6.3.2, $\delta
     =1$]{Bobienski-N:SO(3)}. Taking $b =0$ gives a $3$-dimensional
commutator instead.
\smallbreak

\item $b_{24} = a  \neq 0, b_{45} = b, b_{48} = c, b^2 + c^2
\neq 0,    b_{6} = -2 \frac{b^2+c^2}{a},b_{9} =
- \frac{b^2+c^2+a^2}{a}$, \\[3pt]
$(-2 \frac{b^2+c^2}{a} e^{24} - \frac{b^2+c^2 +
a^2}{a} e^{35},
       -2 b e^{24}, a e^{15} - b e^{25} + c
e^{45},   2 c e^{24},   -a e^{13} + b e^{23}  +
ce^{34})$
$$
      T = 2 \tfrac{b^2+c^2}{a} e^{124}+
\tfrac{b^2+c^2 + a^2}{a} e^{135},\q   dT = 4 \tfrac{(b^2+c^2) (a^2 +
b^2+c^2 )}{a^2}
e^{2345}.
$$
This has $4$-dimensional commutator such that 
${\mathfrak l}^2\iso \so(3)$, and
    curvature $k(\p)=-\bigl(\frac{b^2+c^2}{a}+1\bigr)$.
It corresponds to \cite[Example
      6.3.2]{Bobienski-N:SO(3)}  with $\delta =0$.

\end{enumerate}

\noindent
The first two Lie algebras are $3$-step solvable
with $3$-dimensional commutator and isomorphic.
The third one  is perfect, i.e.  ${\mathfrak l}^1 = {\mathfrak l}$.
For the Lie algebras (1)--(3), and for the family (4) with $a =0$,  ${\mathfrak p}$
is Abelian. The algebra $\mathfrak
h$ is  always $\s_3$, except for (3) where 
${\mathfrak h} \cong {\mathfrak {so}} (3)$.  The first three 
instances fibre onto a Euclidean surface $L/H$,
for the remaining ones $k(\p)$ is negative.
Suitable normalisation of (3) and (5)
give precisely $\gh{l_2},\gh{l_3}$ of Theorem
\ref{thm:main}. 
\end{proof}

The torsion has type $\Lambda^2_3$ for $(4)$ with
$-\sqrt{3}ac=b^2+c^2$, and $(5)$ with $a^2=3(b^2+c^2)$. Up to
$\SO(3)$ equivalence these appear, though in disguise, in
\cite[Theorem 6.7]{Bobienski-N:SO(3)}.
\bigbreak
Let $\widetilde{\Ric}, \Ric$ indicate the Ricci tensors of
the characteristic and Levi-Civita connections. Although in general
there is no reason for the former to be symmetric, this is what
happens at present 
%
\begin{proposition}
If $(\l,\h)$ is a symmetric pair,
the torsion is always coclosed.
\end{proposition}
\begin{proof}
One checks without effort that all Lie algebras in this section have a
symmetric Ricci curvature $\widetilde{\Ric}$. The general formula
for a metric connection with skew torsion
$$
\widetilde{\Ric} (X,Y)=
{\Ric} (X,Y)+\tfrac 12 \hodge d\hodge T(X,Y)+\tfrac 14
\tx\sum_{i,j}^{1,5}
g(T(X,e_i),T(Y,e_j))
$$
allows to conclude.
\end{proof}
\noindent
The same reasoning holds in the non-symmetric case. Since the 
characteristic connection will be identically zero there, every 
$\SO(3)$-curvature tensor vanishes (theorem \ref{thm:non-symm-T=0}),
turning $\hodge T$ into a closed form.
\bigbreak
Lemmas \ref{lemma:symm-non-closed}, \ref{lemma:symm-closed} can
be modified and proved differently, for it is actually possible to spot the closure of $dT$ a
priori. A lengthy computation shows that
$\gh{S}_{x,y,z}(\tilde\na_xT)(y,z,w)=(\tilde\na_wT)(x,y,z)$
everywhere, so from \eqref{eq:nablaT-dT} one has $dT=\tsum
(e_i\hook T)^2$. Now $T$ can be either decomposable as $h\w e^{ij}$, for
some $h\in\h$. The contraction with any basis element is itself 
decomposable and simple, so $\sigma_T=0$ eventually. Therefore the
torsion must be strong by corollary \ref{cor:IP}.
Alternatively, $T$ is proportional to $e^1\w(e^{24}+e^{35})$, yielding a non-zero differential.

\section{Parallel forms}

\subsection{Parallel torsion}

\noindent
Whenever in presence of a characteristic connection, it is relevant to
consider whether it annihilates the
torsion $3$-form. Evidence of this very
restrictive condition can be found in
\cite{Cleyton-S:intrinsic, Alexandrov-FS:aH-revisited}. In the former 
Cleyton and Swann prove 
that parallel torsion implies $d\hodge T=0$, which seems to pervade
our classification. Precisely
\begin{proposition}
Let $(\l,\h)$ be a symmetric pair. Then $\l$ admits
(non-zero) parallel characteristic torsion iff 
it is isomorphic to 
$\so(3)\oplus\R^2$.
\end{proposition}
\begin{proof}
As for lemma \ref{lemma:symm-closed},
when $T$ is closed (and coclosed) the Lie algebras (2) and (4) 
fulfill $\tilde\nabla T=0$. An obvious coordinate change proves them
both Abelian extensions of $\so(3)$. If $dT\neq 0$ instead 
(see lemma \ref{lemma:symm-non-closed}), the torsion is never parallel.
\end{proof}
\bigbreak
\subsection{Reduced holonomy}

\noindent
The holonomy of the characteristic connection will reduce to a subgroup 
only
in presence of a parallel vector, and because of lack of space inside
$\SO(3)$ there is one non-trivial case, that of the circle

\begin{center}
$\Hol(\tilde\na)= \U(1)\ \iff \ \exists\xi\in \l\,:\,\tilde\na\xi=0.$
\end{center}

\noindent
The next result concerns groups $L$ other than the $5$-torus:
\begin{proposition}
Let $(\l,\h)$ be a symmetric pair whose Lie group $L$ has 
characteristic holonomy $\U(1)\sset \SO(3)$ by way of a
$\tilde\na$-parallel vector field $\xi$. Then $\l$ is isomorphic to

1) $\so(3)\oplus\R^2$\smallbreak

2) $(0,e^{14}, e^{15}, e^{12},e^{13})$,\smallbreak

3) $\bigl(-b_{48}(e^{24} + e^{13})+b_9b_{24}e^{35}, 0, e^{15}, e^{24},
e^{31}\bigr)$ 
with $b_{24}=-2(b_{45}^2+b_{48}^2)/b_6$\smallbreak

4) $(0,0,0,0,e^{15}+e^{23})$\smallbreak

5) $\s_3\oplus\R^2$\smallbreak

6) $(0,0,e^{15},e^{53}, e^{34})$.\smallbreak

All others have full holonomy.
\end{proposition}
\begin{proof}
First of all a parallel-vector-to-be $\xi$ must belong to $\h$: 
the
coefficients' relations imply that components in $e_3,e_5$ 
cannot appear. 
Disregarding situations where the connection itself is zero, imposing
$\tilde\na\xi=0$ yields either that 

$\cdot$  $\p$-derivatives are zero (precisely when $\xi=e_1$), or 

$\cdot$ $\tilde\na_\h\equiv 0$ (i.e.~$\xi$ has also
$e_2$-,$e_4$-components).
\smallbreak

\noindent
 a) $e_1$ is parallel for three types of algebras:
\[
\ba{l}
(0,0,b_{24}e^{15}+b_{48}e^{45},-b_{48}e^{35}, -b_{24}e^{13}+b_{48}e^{34})\\[4pt]
\bigl(b_9e^{35}, (b_{19}-b_{45})e^{24}+b_{19}e^{35}, b_{24}e^{15}-b_{45}e^{25}+b_{48}e^{45}, 
-b_{48}e^{35}, -b_{24}e^{13}+b_{45}e^{23}+b_{48}e^{34}\bigr) \\[4pt]
(b_6e^{24},-b_6e^{14},0,b_6e^{12},0)\\[4pt]
(0,b_{13}e^{14}, b_{24}e^{15}, -b_{13}e^{12},-b_{24}e^{13}),\\[4pt]
(b_6e^{24} + b_9e^{35}, -2b_{45}e^{24}, b_{24}e^{15}-b_{45}e^{25}+
b_{48}e^{45}, 2b_{48}e^{24}, b_{24}e^{13}+b_{45}e^{23}+b_{48}e^{34}) 
\ea
\]
with $b_{24}=-2(b_{45}^2+b_{48}^2)/b_6$. The first three are isomorphic.
\smallbreak

\noindent
b) When the line $\xi$ intersects the plane $\lan e_2,e_4\ran$, 
the candidates for covariantly
constant fields have the form 
\[
xe_1+ye_2+te_4,
\]
subjected in particular to $t^2+y^2=3x^2\neq 0$ 
(we omit other relations imposed on \eqref{eq:char-connection}). Apart from
the Abelian algebra  $\gh{a}_5$, here are the solutions:
\[
\ba{ll}
(b_9e^{35},0,b_9e^{51},\tfrac {\sqrt{3}}3 b_9e^{35}, 
b_9e^{13}-\tfrac {\sqrt{3}}3 b_9e^{34}),
& (b_6e^{24},0,b_6e^{41},0, b_6e^{12}),\\[4pt]
\bigl(0,0,b_{24}(e^{15}+\sqrt{3}e^{45}),0,0\bigr),
& \bigl(0,0,0,0,b_{45}(\sqrt{3}e^{15}+e^{23}+e^{45})\bigr),\\[4pt]
(0,b_{19}e^{35},b_{19}e^{52}+\sqrt{3}b_{39}e^{45}, b_{39}e^{35}, 
b_{19}e^{23}+b_{39}e^{34}). 
\ea
\]
Recognising the isomorphisms is easy now.
\end{proof}



\section{ $({\mathfrak l}, {\mathfrak h})$ non-symmetric pair}
\label{sec:non-symmetric}

\noindent
All examples in this section  will have $\tilde \nabla=0$.  They 
have in particular flat characteristic  connection, so
their structure coincides with
 \cite[eqns (6.10)]{Bobienski-N:SO(3)}. 
We will also determine which 
ones have  torsion of type $\Lambda^2_7$.
\bigbreak
Remember that the $SO(3)$-connection  $\tilde\nabla$ satisfies
\eqref{SO(3)conditions} if and only if\smallbreak
\begin{center}
$b_{43} = b_{10}, b_{23} = - b_8,  b_{15}= b_{17} = b_{29} = 
b_{40} = b_{49} =0$.
\end{center}
\smallbreak
\noindent
One can assume at least one of $b_7, b_8,
b_{10}, b_{12},b_{35}, b_{37}$ to be non-zero to prevent 
${\gh p}$ from being $\ad({\gh h})$-invariant and ${\gh h}$  a
Lie subalgebra of ${\mathfrak l}$.
The structure equations then are
\bdm
\left\{
\begin{array}{lcl}
  d e^1 &= &- b_{12} e^{23} + b_6 e^{24} + b_7 e^{25} + b_8 e^{34} +
  b_9 e^{35} + b_{10} e^{45},\\[3pt]
  d e^2 &= &b_{12} e^{13} + b_{13} e^{14} - b_7 e^{15} + 2( b_{50} -
  b_{45} +  b_{19})e^{24} + b_{35} e^{34} + b_{19}
  e^{35} - b_{37}e^{45},\\[3pt]
  d e^3 &=& - b_{12} e^{12} +{\sqrt{3}} b_{28} e^{13} - b_8 e^{14} + b_{24}
  e^{15} + \frac{{\sqrt{3}}}{3} (b_9 + b_{24}-\frac{1}{2}b_{13}
  -\frac{1}{2} b_6)e^{23} - b_{35} e^{24}+\\[3pt]
  && (b_{50} + b_{28} - b_{45}) e^{25} + b_{28} e^{34} + (
  \frac{{\sqrt{3}}}{3} b_9 + \frac{{\sqrt{3}}}{3}
  b_{24} -
  \frac{{\sqrt{3}}} {6} b_{13} - \frac{{\sqrt{3}}}{6} b_6 + 
  b_{48} - b_{47}) e^{45},\\[3pt]
  d e^4 &=& - b_{13} e^{12} + b_8 e^{13} - b_{10} e^{15} + b_{35}
  e^{23} + 2( b_{39} +  b_{48} -  b_{47}) e^{24} + b_{37} e^{25} + b_{39}
  e^{35},\\[3pt]
  d e^5 &=& b_7 e^{12} + ( b_9
  + {\sqrt{3}} b_{47}- \frac{1}{2} b_{13} - \frac{1}{2} b_6) e^{13} + 
  b_{10} e^{14} + {\sqrt{3}} b_{50} e^{15} +
  b_{45} e^{23}  - b_{37} e^{24} + \\[3pt]
  &&  b_{47} e^{25}+ b_{48} e^{34} + b_{50} e^{45},
\end{array}
\right.
\edm
\noindent
and the torsion reads
\begin{equation}
\label{eq:ns-T}
\begin{array} {lcl}
T & = & b_{12} e^{123} - b_6 e^{124} - b_7 e^{125} - b_8 e^{134} + (b_{24}
- \tfrac{1}{2} b_{13}  - \tfrac{1}{2} b_6 - {\sqrt{3}} b_{47} )
e^{135} -\\[3pt]
 & & b_{10} e^{145}- b_{35} e^{234}+ (b_{28}-b_{50}-b_{19} e^{235}+ b_{37} e^{245}+\\[3pt]
 & &  \tfrac{{\sqrt{3}}}{3}( \tfrac 12 b_{13}+ \tfrac 12 b_6 -b_9 - b_{24}+b_{39}
-b_{47}) e^{345}.
\end{array}
\end{equation}

\begin{theorem} 
\label{thm:non-symm-T=0}
If  $(\l, \h)$ is non-symmetric,
then $\tilde \nabla =0$ and $T$ is harmonic
\begin{center}
$dT=0,\qq d\hodge T=0$.
\end{center}
 In addition, $\l$ has
$3$-dimensional commutator ${\l}^1 = {\l}^2$.
\end{theorem}

\begin{proof}
By imposing $d^2 =0$ one gets the extra conditions:
\begin{center}
$b_{50} = b_{47} = b_{28} = 0, b_{48} = - b_{39}, b_{13} = - b_6, b_{24} =
- b_9, b_{45} = b_{19}$,
\end{center}
and then $\tilde \nabla =0$ by (\ref{eq:gamma-connection}). 
Thus the above equations reduce to $(6.10)$ of 
\cite{Bobienski-N:SO(3)} with
$$
t_1 = - b_{12}, t_2 = b_6, t_3 = b_7, t_4 = b_8, t_5 = b_9, t_6 =
b_{10}, 
t_7 = b_{35}, t_8 = b_{19}, t_9 =0, t_{10} = - b_{39}
$$
subject to the constraints
$$
\left \{ \begin{array} {l} -b_{10} b_{19} +b_{7} b_{39}- b_9 b_{37} =0,\\
  b_{7} b_{35} - b_{6} b_{19} + b_{12} b_{37} =0,\\
  -b_{12} b_{39} 
+b_{8} b_{19} -b_{9} b_{35} =0, \\
-b_{6} b_{39}+b_{10} b_{35} + b_8 b_{37} =0,\\
 -b_{7} b_{8} +b_{12}Êb_{10} +b_{6} b_{9} =0.
 \end{array} 
\right.
$$
This system implies $dT = d (\hodge T) =0$. 
This can be alternatively seen by computing the form $\sigma_T$ of
\eqref{eq:sigmaT}, as done previously. 
\bigbreak
A standard computation yields
\begin{enumerate}

\item  $b_6 = - \frac{(-b_7 b_8 +b_{12} b_{10})} {b_9}, b_{37} =  
\frac{(b_7 b_{39} -b_{19} b_{10})} {b_9},
b_{35} = - \frac{(-b_{19} b_8 +b_{12} b_{39})} {b_9}, b_9 \neq 0,$

\item $b_9 =0, b_{35} = - \frac{(-b_6 b_{19} +b_{12} b_{37})} {b_7},
  b_8 =  
\frac{b_{12} b_{10}} {b_7},
b_{39} =  \frac{b_{19} b_{10}} {b_7}, b_7 \neq 0,$

\item  $b_7 =0, b_9 =0, b_{12} =0, b_{19} =0, b_{35} =
\frac{ (-b_8 b_{37} +b_6 b_{39} )} {b_{10}}, b_{10} \neq 0,$ 

\item $b_7 =0, b_{39} =  \frac{b_{12} b_{37}} {b_6}, b_9 = b_{10} =0,
  b_{39} =  
\frac{b_{8} b_{37}} {b_6}, b_6 \neq 0,$

\item $b_{19} =  \frac{b_{12} b_{39}} {b_8}, b_6 = b_7 = b_9 = 
b_{10} = b_{37} =0, b_8 \neq 0,$

\item $b_6 = b_7 =  b_8 =b_9 = b_{10} = b_{12} =0,$

\item $b_6 =b_7 =  b_8 =b_9 = b_{10} = b_{39} =0,$

\end{enumerate}

\noindent
whose corresponding Lie algebras have all  $3$-dimensional
commutator 
${\mathfrak l}^1=[{\mathfrak l}^1, {\mathfrak l}^1].$ 
Tedious details have been left out.
\end{proof}
\bigbreak
Since $\tilde \nabla$ is flat, \cite[proposition
  6.6]{Bobienski-N:SO(3)} has  $\l$ isomorphic to
${\mathfrak {so}}(3) \oplus \R^2$ provided that 
$b_{39}\neq 0$ in \eqref{eq:ns-T}. 
We can generalise this result  regardless of 
coefficients

\begin{theorem}
$\tilde \nabla =0$ on $\l$ if and only if $\l$ is isomorphic to $\so(3)\oplus\R^2$.
\end{theorem}
\begin{proof}
The connection is zero if and only if $g([X,Y], Z)$ is a
totally skew. In other words if 
$$
g (\ad_X Y, Z) = - g (Y, \ad_X Z),
$$
for any $X, Y, Z \in \l$. Then $\l$ is compact and Weyl's theorem says that
 ${\l} ={\l}^1 \oplus {\mathfrak z}$, where ${\mathfrak z}$ is the
center of $\l$. By theorem \ref{thm:non-symm-T=0}
 the commutator ${\l}^1$ is a
$3$-dimensional centerless compact Lie algebra, whence 
semisimple. Thus it isomorphic to ${\mathfrak {so}}(3)$.
\end{proof}

\begin{rmk1}
The characteristic connection being zero is reminiscent of compact
semisimple Lie groups of even dimension equipped with the complex
structure of Samelson and metric equal to the negified Killing
form. They have in fact zero Bismut connection and the corresponding
torsion is always harmonic, exactly as in \ref{thm:non-symm-T=0}, 
see  \cite{Spindel-STvP:complex}.
\end{rmk1}
\begin{rmk2}
Note incidentally that if $\tilde\na\equiv 0$ then $\l=\h\oplus\p$
defines a naturally reductive space, for both the characteristic
curvature and the torsion are obviously parallel \cite{Tricerri-V:naturally-reductive}.
\end{rmk2}
\bigbreak
As equations \eqref{eq:ns-T} translate now in to 
$$
\sqrt{3} b_{35}=  b_{10} - b_{12}, \quad 
b_8 + b_7 = \sqrt{3} b_{37}, \quad b_6 = - 2 b_9,
$$

\begin{corollary}
If  $(\l,\h)$ is non-symmetric,  the torsion $T$ is of pure type 
$\Lambda^2_7$ when:
\begin{enumerate}

\item $b_6 = - 2 b_9,  b_{12} = - \sqrt{3}b_{35} +b_{10}, 
b_{39} = \frac{\sqrt{3}}{6 b_9} (2 b_9^2- b_{10}^2-b_8^2), 
b_{37} = -\frac {{\sqrt{3}}} {3 b_8} (- b_8^2+
2 b_9^2-3 b_{10}^2+{\sqrt{3}} b_{10} b_{35}), 
b_{19} = \frac{8\sqrt{3}}{3}(- b_{10} b_{8} ^2+
2 b_{10} b_9^2- b_{10} ^3+ {\sqrt{3}} b_{10}^2 b_{35} +{\sqrt{3}} b_{8}^2 b_{35}), 
b_7 = - \frac{1}{b_8} (2 b_{9}^2-b_{10}^2 + \sqrt{3} b_{10} b_{35})$

\item $b_6 = - 2 b_9, b_8 =0, b_{19} = - \frac{b_{10} b_{37}} {2 b_9}, 
b_{39} = \frac{\sqrt{3}}{6 b_9} (2b_9^2-b_{10}^2), 
b_{35} = - \frac{\sqrt{3}}{3 b_{10}} (2b_9^2-b_{10}^2), 
b_{12} = \frac{2 b_9^2}{b_{10}}, b_7 = \sqrt{3} b_{37}$

\item $b_8 =  b_{35} = b_{39} =0, b_{10}^2 = 2 b_9^2,  
b_6 = - 2 b_{9}, b_{12} = b_{10}, b_7 = \sqrt{3} b_{37}, b_{37}^2 = 2 b_{19}^2$

\item $b_6 =  b_8 = b_9 = b_{10} = b_{39} =0, 
b_{12} = - \sqrt{3} b_{35}, b_7 = \sqrt{3} b_{37}$

\item $b_6 = b_7 = b_8 =  b_9 = b_{10}  = b_{12} = b_{35} = b_{37} =0$.\qed 

\end{enumerate}
\end{corollary}

\section{Example of pure type $\Lambda^2_7$ with  non-closed torsion}
\label{sec:example} 

\noindent
The majority of Lie algebras found are strongly \NI, 
in the sense that $T$ is closed. 
It is tempting to think this feature could be proven in general, but
this is not the case.
The disproving example is constructed using a fibration similar to
\eqref{eq:fibration} and will be  required to have type  $\Lambda^2_7$.
\smallbreak
We consider Lie
algebras $\l$ whose Levi-Civita connection satisfies
\begin{equation} \label{inverseLeviCivita}
\begin{array}{l}
\nabla_X {\gh h} \subseteq {\gh p}, \quad \nabla_X {\gh p}
\subseteq {\gh h}, \quad \forall X \in {\gh h}\\
\nabla_Y {\gh h} \subseteq {\gh h},  \quad \nabla_Y {\gh p}
\subseteq {\gh p}, \quad \forall Y \in {\gh p}.
\end{array}
\end{equation}
These are in some sense \lq dual\rq\  to those of proposition
\ref{prop:LC conditions}. The constraints \eqref{inverseLeviCivita} are equivalent 
to $({\l}, {\p})$  being a symmetric pair this time, i.~e.
$$
[\p, \p] \subseteq \p, \quad [\h, \p] \subseteq \h, \quad [\h, \h] \subseteq \p.
$$
%
The $1$-connected Lie subgroup $P\sset L$ with Lie algebra 
$\p$ is the typical fibre of the Riemannian submersion
$$
L \longrightarrow L /P
$$
that renders $L/P$ a 
$3$-dimensional locally symmetric space. With the aid of O'Neill
formulas the Ricci tensor can be expressed via 
$$
\textsl{Ric}(X, X) = - ([X, [X, e_1]_{\mathfrak p} ]_{\mathfrak h},
e_1) - ([X, [X, e_2]_{\mathfrak p} ]_{\mathfrak h},
e_2) - ([X, [X, e_4]_{\mathfrak p} ]_{\mathfrak h},
e_4),
$$
for any $X \in {\mathfrak h}$, and yields the scalar curvature in terms of the $b_i$'s.
\bigbreak
%
%
%
The Lie algebra 
 $$
\l=( - \tfrac 34  e^{15} + \tfrac {3\sqrt{3}}4  e^{23} - 
\tfrac {\sqrt{3}}4   e^{45},  e^{25},  
\sqrt{3}     e^{12}  -  e^{24} -  e^{35}, 
- \tfrac 54  \sqrt{3}  e^{15} - \tfrac 94   e^{23} - \tfrac 54  e^{45}, 0)
 $$
is a particular solution of the system \eqref{eq:pure} coupled with
the Jacobi identity. It has highest step length (four) and commutator
${\l}^1 = \text{span}\{e_1, \ldots, e_4\}$. 
%
Its torsion 
\[
T =  \tfrac {\sqrt{3}}4  e^{123} - \sqrt{3}    e^{145} - \tfrac 34   e^{234}
\]
is not strong 
\[
dT =  - \tfrac 32 ( e^{2345} +   \sqrt{3}  e^{1235}),\q d\hodge T=0.
\]
The Ricci curvature 
\[
\textsl{Ric}(L)=
 \left(
 \begin{matrix}
   - \tfrac{141} {32}&0&0&- \frac{99} {32} \sqrt{3}&0\\
   0&-\frac{27} {8}&0&0&0\\
   0&0&-\frac{27} {8}&0&0\\
   - \frac{99} {32} \sqrt{3}&0&0&- \frac{99} {32} \sqrt{3}&0\\
0&0&0&0&- \frac{15} {2}
 \end{matrix}
 \right)
 \]
is negative, and $L/P$ is scalar-flat.

\section{Higher dimensional geometry}\label{sec:8}

\noindent
Although $\SO(3)$ geometry is quite intriguing on its own, it is the
links with other $G$-structures that make it really valuable.
From example \ref{ex:hypo} in fact, one can build $6$-manifolds with
holonomy $\SU(3)$, a feature shared by probably many other instances
of the same kin.
 
A direct jump to dimension $7$ is possible, and twofold alluring: 
first because there exist non-integrable
CR-structures arising from $5$-manifolds of pure type
$\Lambda^2_7$. Secondly, dimension seven is also inhabited by $\G$
metrics, and it is all the more natural to consider $2$-sphere bundles 
over $L$ related to a twistor theory of sorts \cite{Bobienski-N:SO(3)}.
\smallbreak

We wish to concentrate on dimension eight now though, motivated by 
\cite{Nurowski:distinguished-dimensions}. There  the author
investigates $G$-structures related to
a series of simple Lie groups 
$\SO(3),\ \SU(3),\ldots$ acting in dimensions
$5,8,\ldots$, describes the corresponding \NI\ conditions and discusses the
existence of a characteristic connection.  
We thus consider on the $8$-dimensional product $L \times K$, with $K = \R^3$ or
$\SO(3)$, the following
left-invariant metric
$$
\tilde g = g +  (e^6)^2 + (e^7)^2 + (e^8)^2,
$$
where $g$ is the left-invariant metric on $L$ of page \pageref{eq:rep}
and $\{e_6, e_7, e_8\}$
is a basis of the Lie algebra $\mathfrak k$ of $K$.
When $K = \SO(3)$ it is the standard basis
$ [e_6, e_7] = e_8, [e_7, e_8] = e_6, [e_8, e_6] = e_7.
$
The vector space $\R^8$ is identified with the set of traceless skew-Hermitian
$3\times 3$ matrices by way of
\bdm
\tilde X = (x_1, \ldots, x_8) \leftrightarrow
\begin{pmatrix}
x_1 - \sqrt{3} x_4 & \sqrt{3} (x_2 + i x_8) & \sqrt{3}(x_3 + i x_7)\\
\sqrt{3} (x_2 - i x_8)& x_1 + \sqrt{3} x_4 & \sqrt{3} (x_5 + ix_6)\\
\sqrt{3} (x_3 - i x_7) & \sqrt{3} (x_5 - i x_6) & - 2 x_1
       \end{pmatrix}.
\edm
The irreducible representation of $\SU(3)$ on $\R^8$ is given -- as
for $\SO(3)$ --  by
$\tilde \rho (h) \tilde X = h \tilde X h^{-1}, h \in \SU(3)$. An
$\SU(3)$ structure on $(L \times K, \tilde g)$ is
defined  \cite{Nurowski:distinguished-dimensions} by  an element
$\tilde{\mathbb T} \in \otimes^3 \R^8$
satisfying similar relations to those of ${\mathbb T}$
$$
\begin{array} {ll}
\tilde{\mathbb T} (\tilde X, \tilde X, \tilde X)  & =
{\mathbb T} (X, X, X) + \tfrac 32 x_1 (x_6^2 + x_7^2 -2 x_8^2) -
\tfrac {3\sqrt{3}}{2}  x_4 (-x_7^2+ x_6^2)\\[3pt]
{}& \phantom{=}+ 3 \sqrt{3}
(x_3 x_6 x_8 - x_5 x_6 x_7 +  x_2 x_7 x_6),\\[3pt]
& =   \tfrac 12 \det\tilde X
\end{array}
$$
with $X = (x_1, \ldots, x_5)$. Setting $x_6,x_7,x_8$ to zero induces
the same structure on $\R^5\sset\R^8$ as \eqref{eq:rep}, because 
the position of $\sqrt{3}$ is only
cosmetic. 
Yet in stark contrast to dimension five \cite{Bobienski-N:SO(3)},
the reduction is equivalently determined by the differential $3$-form
\[
\psi = E_1 \wedge e^6 + E_2 \wedge e^7 + E_3 \wedge e^8 + e^{678},
\]
where the $E_j$'s come from \eqref{eq:E-forms}. The isotropy of $\psi$
is $\SU(3)/\bb{Z}_3$ embedded in  $\GL(8,\R)$. We refer the
reader to \cite{Hitchin:forms} where the action of $\SU(3)$ on $\R^8$ 
was first examined in detail.
All details in the same flavour can be found in \cite{Witt:triality}.
As in the lower dimension, a nearly
integrable $\SU(3)$ structure is defined in terms of a 
 symmetric Killing tensor
$\tilde{\mathbb T}$. 
This only implies the existence of a connection with totally skew
torsion, but is no longer equivalent to it.
We shall see under which circumstances one can get hold of a nearly
integrable $\SU(3)$ structure  on $L\times K$, and 
that the induced characteristic connection is the zero connection,
making $\psi$ obviously parallel. 

The explicit components of $\tilde{\mathbb T}=\sum_{i,j,k=1}^8\tilde
t_{ijk}\,dx_i\otimes dx_j\otimes dx_k$ are:
\begin{equation}
    \label{eq:tt-components}
\ba{c}
\tilde t_{ijk}=t_{ijk}, \q i,j,k\leq 5;\\[3pt]
\tilde t_{166} =-\tfrac 12=\tilde t_{177},\q \tilde t_{188}=1,\q 
\tilde t_{466}=-\tilde t_{477}=\tilde t_{267}=-\tilde t_{368}=
\tilde t_{578}=\tfrac {\sqrt{3}}2.
\ea
\end{equation}
Let $\nabla$ indicate the Levi-Civita connection on $\l + \k$
and also its restriction to $\l$ and $\k$, so 
\bdm
\nabla_{e_6} e_7 = \tfrac 12  e_8=-\nabla_{ e_7} e_6,\q
\nabla_{ e_6} e_8 = -\tfrac 12  e_7=-\nabla_{ e_8} e_6,\q
\nabla_{ e_7} e_8 = \tfrac 12  e_6=-\nabla_{ e_8} e_7
\edm
for $K =\SO(3)$. 
\begin{lemma}
Let $(L^5,g,{\mathbb T})$ be \NI.
Then the $\SU(3)$ structure $(L\times K,\tilde g, {\tilde{\mathbb T}})$,
with $K=\R^3$ or $\SO(3)$ is \NI\ if and only if
\[
\nabla_X X=0\q \text{for all}\ X\in\l.
\]
\end{lemma}
\begin{proof}
Set $\tilde X = X + Y\in\l+\k$. 
By construction $\nabla_X Y = 0$, and both when $\k = \R^3$ and
$\k = \so(3)$, $\nabla_Y Y $ vanishes as well. Thus
\begin{align*}
0&={\tilde{\mathbb T}} (\nabla_{\tilde X} \tilde X, \tilde X, \tilde
X) = {\tilde{\mathbb T}} (\nabla_X X, \tilde X, \tilde X) \\
& =
{\mathbb T}(\nabla_X X, X, X) + 2 {\tilde{\mathbb T}} (\nabla_X X, X, Y) +
{\tilde{\mathbb T}} (\nabla_X X, Y, Y).
\end{align*}
The first term (all $X$'s) vanishes by near integrability of $L$, and the
second too because the components \eqref{eq:tt-components} yield
${\tilde{\mathbb T}} (\l, \l, \k) =0$. By linearity 
${\tilde{\mathbb T}}(\nabla_X X, e_i, e_j) =0,\ i,j=6,7,8$. Using
\eqref{eq:tt-components} once more gives $g(\nabla_X X, \l)=0$.
The argument also works backwards.
\end{proof}

\noindent
The relation of the lemma says that $\nabla_{e_i} e_j +
\nabla_{e_j} e_i = 0$ for all $i,j=1,\ldots, 5$, which implies
\eqref{eq:flatness}. The first string of which characterises
\eqref{SO(3)conditions}, so
\begin{corollary}
The geodesic equation $\nabla_X X=0$ induces the split behaviour of
$\tilde\na$ on $\l$ described by \eqref{SO(3)conditions} and forces
$\tilde \nabla$  to be zero.\qed
\end{corollary}

In this case, the skew-symmetry of $\nabla$ 
can be simply detected by considering the relation 
$2g(\nabla_{e_i}e_j + \nabla_{e_j}e_i, e_k) = 
T(e_i, e_j, e_k) + T(e_j, e_i, e_k)$, $\forall i,j,k$.

\noindent
That said, generating nearly integrable $\SU(3)$ structures on the
$8$-manifold $(L\times K,\tilde g, {\tilde{\mathbb T}})$ is 
possible only by means of one  candidate
\begin{proposition}
Up to isomorphisms, the unique symmetric pair giving rise to 
\NI\ products $L\times \R^3, L\times \SO(3)$ 
is  $\bigl(\l=\so(3)\oplus\R^2,\h=\so(3)\bigr)$.
\end{proposition}
\begin{proof}
In the symmetric case  $\nabla_X X=0$ reduces the structure equations
of $L$ to
\bdm
\left\{
\begin{array} {l}
d\1 = b_6 e^{24} + b_9e^{35}\\
d\2 = - b_6 e^{14} + b_{19} e^{35}\\
d\3 =  e^5\w(b_9 e^{1} + b_{19} e^{2} + b_{39} e^{4})\\
d\4 = b_6 e^{12}+b_{39}e^{35}\\
d\5 =  (b_9 e^{1} + b_{19} e^{2} + b_{39} e^{4})\w e^{3}.\\
\end{array}
\right.
\edm
From the Jacobi identity 
either $b_6 =0$ 
 or $b_9 = b_{19} 
=b_{39}$.
Both structures are isomorphic to $\so(3)\oplus\R^2$
and crop up as II, IV in lemma \ref{lemma:symm-closed}.
\end{proof}

\begin{proposition}
Non-symmetric pairs with $b_{37} =0$ always induce a nearly
integrable $\SU(3)$ structure on the product $L\times K$.
\end{proposition}
\begin{proof}
All algebras of section \ref{sec:non-symmetric} satisfy, by way of
$d^2=0$, the extra requirements
$b_{50} = b_{47} = b_{28} = 0, b_{48} = - b_{39}, b_{13} = - b_6, b_{24} =
- b_9, b_{45} = b_{19}$.
\end{proof}


\vspace{0.6cm}

\end{document}